\newif\iftc 

\if== 
 \documentclass[12pt,a4paper]{article}
 \tcfalse 
\fi

\if!=      
\usepackage[colorlinks, bookmarks, bookmarksopen, bookmarksnumbered,
 bookmarksopenlevel=3,pdfstartview=FitH]{hyperref}[1999/08/31]
\usepackage{ttftimes}
\else
\fi

\usepackage{amssymb}
\usepackage{amsfonts}
\usepackage{latexsym}
\usepackage{theorem}

\usepackage{amscd}
\usepackage{amsfonts}
\usepackage{stmaryrd}
\usepackage{eufrak}
\usepackage{euscript}
\usepackage{mathrsfs}

\if==                     
\usepackage{cite}
\fi

\newcommand*{\qed}{\penalty\@M\qquad\mbox{$\square$}}  

\theoremheaderfont{\upshape\rmfamily\bfseries}

\theorembodyfont{\slshape}

    \newtheorem{trm}{Theorem}[section]
    \newtheorem{prp}[trm]{Proposition}
    \newtheorem{lem}[trm]{Lemma}

\theorembodyfont{\upshape}

    \newtheorem{dfn}[trm]{Definition}

\newenvironment{theorem}{\trm}{\endtrm}
\newenvironment{lemma}{\lem}{\endlem}
\newenvironment{proposition}{\prp}{\endprp}
\newenvironment{definition}{\dfn}{\enddfn}

\begin{document}


\newcommand{\BF}[1]{{\bf#1}}
\newcommand{\SF}[1]{{\sf#1}}
\newcommand{\RM}[1]{{\rm#1}}

\renewcommand{\SS}{\subseteq}

\newcommand{\NW}[1]{{\em#1}}
\newcommand{\UP}[1]{{\rm#1}}


\newcommand{\lf}[1]{\label{F:#1}}  
\newcommand{\rf}[1]{{\rm(\ref{F:#1})}}  
\newcommand{\ls}[1]{\label{S:#1}}  
\newcommand{\rs}[1]{\ref{S:#1}}  
\newcommand{\lp}[1]{\label{P:#1}}  
\newcommand{\rp}[1]{\ref{P:#1}}  
\newcommand{\ct}[1]{\cite{#1}}          


\newcommand{\Iff}{\;\iff\;}
\newcommand{\iffdef}{\buildrel{\RM{def}}\over\iff}
\newcommand{\Iffdef}{\;\iffdef\;}
\newcommand{\imp}{\;\Longrightarrow\;}
\newcommand{\Imp}{\;\imp\;}
\newcommand{\AND}{\;\;\&\;\;}
\newcommand{\eqdef}{\buildrel{\RM{def}}\over=}
\newcommand{\Eqdef}{\;\eqdef\;}



\newcommand{\cls}[1]{{\BF{#1}}}     

\newcommand{\SLS}{\cls{SLS}}                       
\newcommand{\TLS}{\cls{TLS}}                       
\newcommand{\EUC}{\cls{EUC}}                       
\newcommand{\MVT}{\cls{MVT}}                       
\newcommand{\Set}{\cls{Set}}                       
\newcommand{\Rel}{\cls{Rel}}                       

\newcommand{\UU}{\bf{U}}
\newcommand{\CC}{{\bf{C}}}
\newcommand{\DD}{{\bf{D}}}

\newcommand{\EC}[1]{{\bf#1}}               


\newcommand {\spc}[1]{{\cal #1}}

\newcommand{\A}{\spc{A}}
\newcommand{\B}{\spc{B}}
\newcommand{\C}{\spc{C}}
\newcommand{\R}{\spc{R}}    
\newcommand{\D}{\spc{D}}    
\newcommand{\U}{\spc{U}}    
\newcommand{\T}{{\spc{T}}}  


\newcommand{\LI}{\preccurlyeq}                
\newcommand{\BI}{\succcurlyeq}                
\newcommand{\EI}{\sim}                        

\newcommand{\LS}{\sqsubseteq}                 
\newcommand{\BS}{\sqsupseteq}                 

\newcommand{\LQ}{\ll}                         
\newcommand{\BQ}{\gg}                         
\newcommand{\EQ}{\approx}                     

\newcommand{\GE}{\geqslant}  
\newcommand{\LE}{\leqslant}  

\newcommand{\HP}{\vartriangleright}          
\newcommand{\LP}{\vartriangleleft}           

\newcommand{\cUp}{{\textstyle\bigcup}}


\newcommand{\CE}[1]{\bar{#1}}


\newcommand{\Ob}{\mathop{\RM{Ob}}\nolimits}
\newcommand{\Ar}{\mathop{\RM{Ar}}\nolimits}
\newcommand{\Opt}{\mathop{\RM{Opt}}\nolimits}


\newcommand{\smalltimes}{\mathbin{\mathchoice%
           {\raise.2ex\hbox{$\scriptstyle\times$}}%
           {\raise.2ex\hbox{$\scriptstyle\times$}}%
           {\raise.0ex\hbox{$\scriptscriptstyle\times$}}%
           {\raise.0ex\hbox{$\scriptscriptstyle\times$}}}}

\renewcommand{\(}{\left(\big.}
\renewcommand{\)}{\right)}

\renewcommand{\.}{\circ}                                   
\renewcommand{\*}{\smalltimes}                             

\renewcommand{\:}{\colon}                                  
\renewcommand{\>}{\to}                                     

\renewcommand{\#}{\times}                                  

\newcommand{\q}{\hskip.7em\relax}

\newcommand{\<}[1]{_{{}_{\!\scriptstyle#1\!}}}          

\newcommand{\It}[1]{\mbox{\UP{(#1)}}}

\renewcommand{\=}{=\\&=&}  

\newcommand{\EMPTYSET}{\varnothing}


\newcommand{\IF}{&\mbox{if}&}

\def\CASE#1&#2;#3&#4;{ \left\{\begin{array}{lll}
                         #1,\IF\!\!#2,\\[2mm]#3,\IF\!\!#4
                       \end{array}\right. }

\makeatletter

\begingroup
  \catcode`\|=13
  \gdef\SET{\catcode`\|=13
\def|{\setbox0=\hbox\bgroup$\displaystyle\;\left\vert\;\mathstrut}\SET@}

\gdef\SET@#1{{\left\{               
  \setbox1=\hbox{$\let|=\;#1$}
  #1 \vphantom{{\let|=\;\displaystyle#1}}
  \right.$\egroup\ht0=1.2\ht0\box0
  \right\} }}
\endgroup

\makeatother

\def\TPL#1{\left\langle#1\right\rangle}


\let\UL=\_

\def\_{\nobreak\hskip0pt-\allowbreak\hskip0pt} 

\newcommand {\be}{\begin{equation}}
\newcommand {\ee}{\end{equation}}
\newcommand{\bea}{\begin{eqnarray}}
\newcommand{\eea}{\end{eqnarray}}
\newcommand{\ba}{\begin{array}}
\newcommand{\ea}{\end{array}}
\newcommand{\beq}{\begin{eqnarray*}}
\newcommand{\eeq}{\end{eqnarray*}}

\def\clX{\spc{X}}
\def\clY{\spc{Y}}
\def\clD{\spc{D}}

\newcommand{\1}{{\bf 1}}
\newcommand{\0}{{\bf 0}}
\newcommand{\bI}{{\bf I}}

\newcommand{\ds}{\displaystyle}
\newcommand{\sss}{\scriptscriptstyle}

\newcommand{\ddef}{\stackrel{\rm def} {\ds=}}

\def\notni{\hbox{{\kern0pt\raise.7pt\hbox{${\scriptstyle \times}$}}
{\kern-9.5pt\raise0pt\hbox{$\supset$}}}}

\def\dotvee{\hbox{{\kern0pt\raise3pt\hbox{${\cdot}$}}
{\kern-9.65pt\raise0pt\hbox{$\vee$}}}}

\def\dotcup{\hbox{{\kern0pt\raise3pt\hbox{${\cdot}$}}
{\kern-9.65pt\raise0pt\hbox{$\cup$}}}}

\def\dotwedge{\hbox{{\kern10pt\raise-2pt\hbox{${\cdot}$}}
{\kern-9.65pt\raise0pt\hbox{$\wedge$}}}}

\def\dotcap{\hbox{{\kern10pt\raise-2pt\hbox{${\cdot}$}}
{\kern-9.65pt\raise0pt\hbox{$\cap$}}}}

\def\dotbigvee{\hbox{{\kern10pt\raise4pt\hbox{${\cdot}$}}
{\kern-10.5pt\raise0pt\hbox{$\bigvee$}}}}

\def\dotbigcup{\hbox{{\kern10pt\raise4pt\hbox{${\cdot}$}}
{\kern-10.5pt\raise0pt\hbox{$\bigcup$}}}}

\def\dotbigcap{\hbox{{\kern10pt\raise-2pt\hbox{${\cdot}$}}
{\kern-10.5pt\raise0pt\hbox{$\bigcap$}}}}

\def\dotbigwedge{\hbox{{\kern10pt\raise-2pt\hbox{${\cdot}$}}
{\kern-10.5pt\raise0pt\hbox{$\bigwedge$}}}}

\def\dotar{\hbox{{\kern4pt\raise.0pt\hbox{${\cdot}$}}
{\kern-9.5pt\raise0pt\hbox{$\Rightarrow$}}}}

\def\dotal{\hbox{{\kern0pt\raise5pt\hbox{${\cdot}$}}
{\kern-9.5pt\raise-0.5pt\hbox{$\forall$}}}}

\def\dotex{\hbox{{\kern0pt\raise3pt\hbox{${\cdot}$}}
{\kern-9.5pt\raise-0.5pt\hbox{$\exists$}}}}

\def\dotin{\hbox{{\kern0pt\raise1.8pt\hbox{${\cdot}$}}
{\kern-10.5pt\raise0pt\hbox{$\in$}}}}

\def\dotsubset{\hbox{{\kern0pt\raise.2pt\hbox{${\cdot}$}}
{\kern-10.5pt\raise0pt\hbox{$\subset$}}}}


\newcommand{\Meas}{\cls{Meas}}

\newcommand{\Z}{{\spc{Z}}}  
\renewcommand{\P}{{\spc{P}}}
\newcommand{\Q}{{\spc{Q}}}
\newcommand{\I}{{\spc{I}}}  

\newcommand{\al}{\alpha}
\newcommand{\la}{\lambda}
\newcommand{\rh}{\rho}
\newcommand{\si}{\sigma}
\newcommand{\de}{\delta}
\newcommand{\et}{\eta}
\newcommand{\ga}{\gamma}

\renewcommand{\o}{\otimes}

\newcommand{\ST}{\cls{ST}}
\newcommand{\SLT}{\cls{SLT}}


\def\Put(#1)#2{\put(#1){\raisebox{-.8ex}{\makebox[0pt][c]{$#2$}}}}
\def\Arr(#1)#2{\put(#1){\raisebox{-.8ex}{\makebox[0pt][c]{\small$#2$}}}}
\def\Vec(#1;#2;#3){\put(#1){\vector(#2){#3}}}
\def\Lin(#1;#2;#3){\put(#1){\line(#2){#3}}}

\def\Dia(#1,#2)#3{%
\countdef\y=1\y=#2%
\countdef\dy=2\dy=2
\advance\y by\dy\advance\y by\dy%
\begin{picture}(#1,\number\y)(0,-\number\dy)%
#3%
\end{picture}%
}

\newcommand\BegDiag{$$}
\newcommand\EndDiag{$$}

\newcommand{\LongEq}[2]{%
\iftc
\begin{eqnarray*}
 \lefteqn{#1}\\[1mm]\lefteqn{#2}
\end{eqnarray*}
\else
\[
  #1 \qquad #2
\]
\fi
}


\newcommand{\DefSth}{%
\renewcommand{\t}{\theta}
\renewcommand{\S}{{\mathfrak S}}
\renewcommand{\O}{\Omega}
\renewcommand{\o}{\omega}
\newcommand{\s}{\sigma}
\def\Int_##1{\int\limits_{##1}\!\!}
\newcommand{\MS}[1]{\TPL{\O_{##1},\S_{##1}}}
}

\newcommand{\DefLin}{%
\renewcommand{\S}{\Sigma}
}

\newcommand{\DefMul}{%
\newcommand{\BIW}{\mathrel{\dot\BI}}     
\renewcommand{\P}{{\spc{P}}}
}

\newcommand{\DefFuz}{%
\newcommand {\FMT}{\cls{FMT}}    
\newcommand {\FPT}{\cls{FPT}}    
\newcommand{\x}[1]{\mu_{{}_{\scriptstyle\vphantom\beta##1\!}}}
\renewcommand{\d}[1]{\delta_{##1}}             
\newcommand{\II}{{\bf{I}}}
}

\begin{center}
{\Large\bf Fuzzy Logic, Informativeness and Bayesian Decision-Making Problems}

\bigskip

{Peter V. Golubtsov \\[4mm]
\it Department of Physics, Moscow State Lomonosov University \\
\it 119899, Moscow, Russia \\[4mm]
\it E-mail: P\UL V\UL G@mail.ru
}

\bigskip

Stepan S. Moskaliuk

\medskip
{\it Bogolyubov Institute for Theoretical Physics\\
Metrolohichna Str., 14-b, Kyiv-143, Ukraine, UA-03143\\
e-mail: mss@bitp.kiev.ua}

\end{center}

\bigskip
\centerline{Abstract}

\medskip

 This paper develops a category-theoretic approach to uncertainty, informativeness and 
decision-making problems. It is based on appropriate first order fuzzy 
logic in which not only logical connectives but also quantifiers have 
fuzzy interpretation. It is shown that all fundamental concepts of 
probability and statistics such as joint distribution, conditional 
distribution, etc., have meaningful analogs in new context. This 
approach makes it possible to utilize rich conceptual experience of 
statistics. Connection with underlying fuzzy logic reveals the logical 
semantics for fuzzy decision making.
Decision-making
problems within the framework of IT-categories and generalizes
Bayesian approach to decision-making with a prior information are considered. It leads to fuzzy Bayesian approach in 
decision making and provides methods for construction of optimal 
strategies. 

\newpage

\section{Introduction}

 The theory of fuzzy sets is used more and more widely in the 
description of uncertainty. Indeed, very often some poorly 
formalizable notions or expert knowledge are readily expressed in terms 
of fuzzy sets. In particular, fuzzy sets are extremely convenient in 
descriptions of linguistic uncertainties \ct{Zad}. On the other hand, fuzzy 
notions themselves often admit flexible linguistic interpretations. 
This makes the exploitation of fuzzy sets especially natural and 
illustrative.

 It is well known that fuzzy set theory suggests a wide range of 
specific approaches to decision problems in a fuzzy environment \ct{Dub:FSI, Dub:FSII}. 
However, it still cannot easily compete with probability theory in 
dealing with uncertainty. Certainly, mathematical statistics (more 
 precisely, statistical decision theory) has accumulated a rich 
conceptual experience.

 It was shown in \ct{Gol:FuzDec, Gol:MSI} that the basic constructions and propositions of 
probability theory and statistics playing the fundamental role in 
decision-making problems have meaningful counterparts in fuzzy 
theories. It makes it possible to use the methodology of statistical 
 decision-making in the fuzzy context. In particular, the fuzzy variant 
of the Bayes principle derived in this paper plays the same role in 
fuzzy decision-making problems as its probabilistic prototype in the 
theory of statistical games \ct{Black:Gir}.

 In contrast to the approach of \ct{Gol:FuzDec} where all the notions were 
introduced ``operationally'' in order to be more close to the similar 
notions in statistics, in this paper all the notions
 are introduced ``logically,'' i.e., by the corresponding formulas in an 
appropriate fuzzy logic of first order. In this logic not only logical 
connectives but also quantifiers have fuzzy interpretation.

 Connection with underlying fuzzy logic provides an interesting logical 
semantics of fuzzy decision making. Indeed in this approach a priory 
information may be represented by a fuzzy predicate and an 
experiment -- by certain fuzzy relation. The loss function is replaced 
 by fuzzy relation ``good decision for,'' mathematical expectation 
operator -- by fuzzy universal quantifier, etc. As a result the notion 
``good decision strategy'' is expressed by a first order formula in this 
logic.

It is convenient to consider different systems that take place in
information acquiring and processing  as particular cases of
so-called \NW{information transformers} (ITs). Besides, it is
useful to work with families of ITs in which certain operations,
e.g., \NW{sequential} and \NW{parallel compositions} are defined.

It was noticed fairly long ago~%
\ct{Sackst:StatEq, MorseSackst:StatIso, Chentsov,
Chentsov:CatMatStat, Gol:MSAlg}, that the adequate algebraic
structure for describing information transformers (initially for
the study of statistical experiments) is
the structure of \NW{category}~%
\ct{MacLane, Herrlich:Strecker, Arrows, Goldblatt}.

Analysis of general properties for the classes of linear,
multivalued, and fuzzy
information transformers, studied in{}~%
 \ct{Gol:MSAlg, Gol:InfCatLinSys, Gol:RelInf, GolFil:MVMS,
 Gol:MV_Inf, Gol:FuzDec, Gol:FuzLog, Gol:CatInfTran},
allowed to extract general features shared by all these classes.
Namely, each of these classes can be considered as a family of
morphisms in an appropriate category, where the composition of
information transformers corresponds to their ``consecutive
application.'' Each category of ITs (or IT-category) contains a
subcategory (of so called, \NW{deterministic} ITs) that has products.
Moreover, the operation of morphism product is extended in a
``coherent way'' to the whole category of ITs.

The work~%
\ct{Gol:Ax-IT,Gol:MS}
undertook an attempt to
formulate a set of ``elementary'' axioms for a category of ITs,
which would be sufficient for an abstract expression of the basic
concepts of the theory of information transformers and for study
of informativeness, decision problems, etc. This paper proposes
another, significantly more compact axiomatic for a category of
ITs. According to this axiomatic a category of ITs is defined in
effect as a \NW{monoidal} category~%
\ct{MacLane,Arrows},
containing a subcategory (of
\NW{deterministic} ITs) with finite products.

Among the basic concepts connected to information transformers
there is one that plays an important role in the uniform
construction of a wide spectrum of IT-categories --- the concept
of \NW{distribution}. Indeed, fairly often an IT $a\:\A\>\B$ can
be represented by a mapping from $\A$ to the ``space
of distributions'' on $\B$ (see, e.g.,~%
\ct{Gol:MV_Inf,Gol:MSI,Gol:FuzDec}%
). For example, a probabilistic
transition distribution (an IT in the category of stochastic ITs) can be
represented by a certain measurable mapping from $\A$ to the space of
distributions on $\B$. This observation suggests to construct a category of ITs
as a \NW{Kleisli} category~%
\ct{MacLane,Kleisli,BarrWells:TTT},
 arising from the following components: an
obvious category of deterministic ITs; a functor that takes an object $\A$ to
the object of ``distributions'' on $\A$; and a natural transformation of
functors, describing an ``independent product of distributions''.

The approach developed in this paper allows to express easily in
terms of IT-categories such concepts as distribution, joint and
conditional distributions, independence, and others. It is shown
that on the basis of these concepts it is possible to formulate
fairly general statement of decision-making problem with a prior
information, which generalizes the Bayesian approach in the theory
of statistical decisions. Moreover, the \NW{Bayesian} principle,
derived below, like its statistical
prototype{}~%
\ct{Borovkov},
reduces the problem of optimal decision strategy
construction to a significantly simpler problem of finding optimal
decision for a posterior distribution.

Among the most important concepts in categories of ITs  is the
concept of (relative) \NW{informativeness} of information transformers.
There are two different approaches to the concept of
informativeness.

One of these approaches is based on analyzing the ``relative
positions'' of information transformers in the corresponding
mathematical structure. Roughly speaking, one information
transformer is regarded as more informative than another one if
with the aid of an additional information transformer the former
one can be  ``transformed'' to an IT, which is similar to (or more
``accurate'' than) the letter one. In fact, this means that all
the information that can be obtained from the latter information
transformer can be extracted from the former one as well.

The other approach to informativeness is based on treating
information transformers as data sources for decision-making
problems. Here, one information transformer is said to be
semantically more informative than another if it provides better
quality of decision making. Obviously, the notion of semantical
informativeness depends on the class of decision-making problems
under consideration.

In the classical researches of Blackwell{}~%
\ct{Black:Comp, Black:EquComp},
the correspondence between
informativeness (Blackwell sufficiency) and semantical
informativeness (Blackwell informativeness) were investigated in a
statistical context. These studies were extended
by Morse, Sacksteder, and Chentsov{}~%
 \ct{Sackst:StatEq, MorseSackst:StatIso, Chentsov:CatMatStat,
Chentsov}, who applied the category theory techniques to their
studies of statistical systems.

It is interesting, that under very general conditions the
relations of informativeness and semantical informativeness (with
respect to a certain class of decision-making problems) coincide.
Moreover, in some categories of ITs it is possible to point out
{\em one} special decision problem, such that the resulting
semantical informativeness coincides with informativeness.

Analysis of classes of equivalent (with respect to
informativeness) information transformers shows that they form a
partially ordered Abelian monoid with the smallest (also neutral)
and the largest elements.

One of the objectives of this paper is to show that the basic
constructions and propositions of probability theory and
statistics playing the fundamental role in decision-making
problems have meaningful counterparts in terms of IT-categories.
Furthermore, some definitions and propositions (for example, the
notion of conditional distribution and the Bayesian principle) in
terms of IT-categories often have more transparent meanings. This
provides an opportunity to look at the well known results from a
different angle. What is even more significant, it makes it
possible to apply the methodology of statistical decision-making
in an alternative (not probabilistic) context.

Approaches, proposed in this work may provide a background for
construction and study of new classes of ITs, in particular, dynamical
nondeterministic ITs, which may provide an adequate description for
information flows and information interactions evolving in time. Besides, a
uniform approach to problems of information transformations may be useful
for better understanding of information processes that take place in complex
artificial and natural systems.

\section{Categories of information transformers} 
\subsection{Common structure of classes of information transformers}
It is natural to assume that for any information transformer $a$ there are
defined a couple of spaces: $\A$ and $\B$, the space of ``inputs'' (or input
signals) and the space of ``outputs'' (results of measurement, transformation,
processing, etc.). We will say that $a$ ``acts'' from $\A$ to $\B$ and denote
this as $a\:\A\>\B$. It is important to note that typically an information
transformer not only transforms signals, but also introduces some ``noise''. In
this case it is \NW{nondeterministic} and cannot be represented just by a mapping
from $\A$ to $\B$.

It is natural to study information transformers of similar type by
aggregating them into families endowed by a fairly rich algebraic
structure~%
\ct{Gol:MSAlg,Gol:InfCatLinSys}.
 Specifically, it is natural to assume that families of ITs
poses the following properties:

(a)  If $a\:\A\>\B$ and $b\:\B\>\C$ are two ITs, then their
\NW{composition} $b\.a\:\A\>\C$ is defined.

(b)  This operation of composition is \NW{associative}.

(c)  There are certain \NW{neutral} elements in these families,
i.e., ITs that do not introduce any alterations. Namely, for any
space $\B$ there exist a corresponding IT $i\<{\B}\: \B\>\B$ such
that $i\<{\B}\.a = a$ and $b\.i\<{\B} = b$.

Algebraic structures of this type are called \NW{categories}~%
\ct{MacLane,Arrows}.

Furthermore, we will assume, that to every pair of information
transformers, acting from the same space $\D$ to spaces $\A$ and $\B$
respectively, there corresponds a certain IT $a*b$ (called \NW{product} of $a$ and
$b$) from $\D$ to $\A\#\B$. This IT in a certain sense ``represents'' both ITs
$a$ and $b$ simultaneously. Specifically, ITs $a$ and $b$ can be
``extracted'' from $a*b$ by means of projections $\pi\<{\A,\B}$ and $\nu\<{\A,\B}$
from $\A\#\B$ to $\A$ and $\B$, respectively, i.e.,
$\pi\<{\A,\B}\.(a*b) = a$,  $\nu\<{\A,\B}\.(a*b) = b$. Note, that
typically, an IT $c$ such that $\pi\<{\A,\B}\.c = a$,  $\nu\<{\A,\B}\.c = b$ is not
unique, i.e., a category of ITs does not have products (in category-theoretic
sense~%
\ct{MacLane, Herrlich:Strecker, Arrows, Goldblatt}%
). Thus, the notion of a category of ITs demands for an accurate
formalization.

Analysis of classes of information transformers studied in{}~%
 \ct{Gol:MSAlg, Gol:MSI, Gol:InfCatLinSys, Gol:RelInf, GolFil:MVMS,
 Gol:MV_Inf, Gol:FuzDec, Gol:FuzLog, Gol:CatInfTran, Gol:MS}, gives
grounds to consider these classes as categories that satisfy
certain fairly general conditions.

\subsection{Categories: basic concepts}  
Recall that a category (see, for example,{}~%
 \ct{MacLane, Herrlich:Strecker, Arrows, Goldblatt}%
) $\CC$ consists of a class of objects $\Ob(\CC)$, a class of morphisms (or arrows) $\Ar(\CC)$,
and a composition operation $\.$ for morphisms, such that:
\begin{itemize}
\item[\It{a}]
To any morphism $a$ there corresponds a certain pair of objects $\A$ and $\B$ (the source and the
target of $a$), which is denoted  $a\:\A\>\B$.
\item[\It{b}]
To every pair of morphisms $a\:\A\>\B$ and $b\:\B\>\C$ their \NW{composition}
$b\.a\:\A\>\C$ is defined.
\end{itemize}

Moreover, the following axioms hold:

\begin{itemize}
\item[\It{c}]
The composition is \NW{associative}:
\[ c\.(b\.a) = (c\.b)\.a. \]
\item[\It{d}]
To every object $\R$ there corresponds an (\NW{identity}) morphism $i\<\R\:\R\>\R$, so that
\[
  \forall a\:\A\>\B,           \qquad
  a\.i\<\A = a = i\<\B\.a.
\]
\end{itemize}

A morphism $a\:\A\>\B$ is called \NW{isomorphism} if there exists a morphism $b\:\B\>\A$
such that $a\.b=i\<\B$ and $b\.a=i\<\A$. In this case objects $\A$ and $\B$ are called
\NW{isomorphic}.

Morphisms $a\:\D\>\A$ and $b\:\D\>\B$ are called \NW{isomorphic} if there exists am
isomorphism $c\:\A\>\B$ such that $c\.a=b$.

An object $\Z$ is called \NW{terminal} object if for any object
$\A$ there exists a unique morphism from $\A$ to $\Z$, which is
denoted $z\<{\A}\:\A\>\Z$ in what follows.

A category $\DD$ is called a \NW{subcategory} of a category $\CC$ if $\Ob(\DD)\SS\Ob(\CC)$,
$\Ar(\DD)\SS\Ar(\CC)$, and morphism composition in $\DD$ coincide with their composition in
$\CC$.

It is said that a category has (pairwise) products if for every
pair of objects $\A$ and $\B$ there exists their \NW{product},
that is, an object $\A\#\B$ and a pair of morphisms
$\pi\<{\A,\B}\:\A\#\B\>\A$ and $\nu\<{\A,\B}\:\A\#\B\>\B$, called
projections, such that for any object $\D$ and for any pair of
morphisms $a\:\D\>\A$ and $b\:\D\>\B$ there exists a unique
morphism $c\:\D\>\A\#\B$, that yields a commutative diagram:
\BegDiag%
\Dia(60,20){
\Put(30,20){\C}
\Put( 0, 0){\A}
\Put(30, 0){\A\#\B}
\Put(60, 0){\B}
\Vec(23, 0;-1,0;20)\Put(15,4){\pi\<{\A,\B}}
\Vec(37, 0; 1,0;20)\Put(45,4){\nu\<{\A,\B}}
\Vec(27,18;-3,-2;24)\Put(13,12){a}
\Vec(33,18; 3,-2;24)\Put(47,12){b}
\Lin(30,17;0,-1;3)\Lin(30,12;0,-1;3)\Vec(30, 7;0,-1;3)\Put(33,10){c}
}
\EndDiag%
i.e., satisfies the following conditions:
\begin{equation} \lf{1}
  \pi\<{\A,\B}\.c = a, \qquad
  \nu\<{\A,\B}\.c = b.
\end{equation}

We call such morphism $c$ the \NW{product of morphisms} $a$ and
$b$ and denote it $a*b$.

It is easily seen that existence of products in a category implies
the following equality:
\begin{equation} \lf{2}
  (a*b)\.d = (a\.d)*(b\.d).
\end{equation}

In a category with products, for two arbitrary morphisms
$a\:\A\>\C$ and $b\:\B\>\D$ one can define the morphism $a\*b$:
\begin{equation} \lf{3}
  a\*b\:\A\#\B\>\C\#\D, \qquad
  a\*b\Eqdef (a\.\pi\<{\A,\B})*(b\.\nu\<{\A,\B}).
\end{equation}
This definition and{}~\rf{1} obviously imply that the morphism
$c=a\*b$ satisfy the following conditions:
\begin{equation} \lf{4}
  \pi\<{\C,\D}\.c = a\.\pi\<{\A,\B}, \qquad
  \nu\<{\C,\D}\.c = b\.\nu\<{\A,\B}
\end{equation}
Moreover, $c=a\*b$ is the only morphism satisfying conditions{}~\rf{4}.

It is also easily seen that{}~\rf{2} and{}~\rf{3} imply the
following equality:
\begin{equation} \lf{5}
  (a\*b)\.(c*d) = (a\.c)*(b\.d).
\end{equation}

Suppose $\A\#\B$ and $\B\#\A$ are two products of objects $\A$ and $\B$ taken in different
order. By the properties of products, the objects $\A\#\B$ and $\B\#\A$ are isomorphic and the
natural isomorphism is
\begin{equation} \lf{6}
  \si\<{\A,\B}\:\A\#\B\>\B\#\A, \qquad
  \si\<{\A,\B}\Eqdef \nu\<{\A,\B}*\pi\<{\A,\B},
\end{equation}
i.e., a unique morphism that makes the following diagram
commutative:
\BegDiag%
\Dia(60,20){
\Put(30,20){\A\#\B}
\Put( 0, 0){\B}
\Put(30, 0){\B\#\A}
\Put(60, 0){\A}
\Vec(23, 0;-1,0;20)\Put(15,4){\pi\<{\B,\A}}
\Vec(37, 0; 1,0;20)\Put(45,4){\nu\<{\B,\A}}
\Vec(27,18;-3,-2;24)\Put(11,13){\nu\<{\A,\B}}
\Vec(33,18; 3,-2;24)\Put(49,13){\pi\<{\A,\B}}
\Lin(30,17;0,-1;3)\Lin(30,12;0,-1;3)\Vec(30, 7;0,-1;3)\Put(35,10){\si\<{\A,\B}}
}
\EndDiag%

Moreover, for any object $\D$ and for any morphisms $a\:\D\>\A$ and $b\:\D\>\B$, the
morphisms $a*b$ and $b*a$ are isomorphic, that is,
\begin{equation} \lf{7}
  \si\<{\A,\B}\.(a*b) = b*a.
\end{equation}

Similarly, by the properties of products, the objects $(\A\#\B)\#\C$ and $\A\#(\B\#\C)$ are
isomorphic. Let
\[
  \al\<{\A,\B,\C}\: (\A\#\B)\#\C\> \A\#(\B\#\C)
\]
be the corresponding natural isomorphism. Its ``explicit'' form is:
\begin{equation} \lf{8}
\iftc 
  \al\<{\A,\B,\C} \eqdef
  (\pi\<{\A,\B}\.\pi\<{\A\#\B,\C})\!*\!
  \((\nu\<{\A,\B}\.\pi\<{\A\#\B,\C})\!*\!\nu\<{\A\#\B,\C}\)\!.
\else 
\fi
\end{equation}
It can be easily obtained with the following diagram:
\BegDiag%
\Dia(80,40){
\Put(20,40){\A\#\B}
\Put(60,40){(\A\#\B)\#\C}
\Put( 0,20){\A}
\Put(40,20){\B}
\Put(80,20){\C}
\Put(20, 0){\A\#(\B\#\C)}
\Put(60, 0){\B\#\C}
\Vec(47,40;-1, 0;20) \Put(38,38){\pi\<{\A\#\B,\C}}
\Vec(33, 0; 1, 0;20) \Put(42, 4){\nu\<{\A,\B\#\C}}
\Vec(17,37;-1,-1;15) \Put( 6,33){\pi\<{\A,\B}}
\Vec(63,37; 1,-1;15) \Put(76,33){\nu\<{\A\#\B,\C}}
\Vec(17, 3;-1, 1;15) \Put( 4,10){\pi\<{\A,\B\#\C}}
\Vec(63, 3; 1, 1;15) \Put(74,10){\nu\<{\B,\C}}
\Vec(23,37; 1,-1;15) \Put(28,28){\nu\<{\A,\B}}
\Vec(57, 3;-1, 1;15) \Put(52,14){\pi\<{\B,\C}}
%
\iftc
\Lin(57,37;-1,-1;5)\Lin(51,31;-1,-1;5)
\Lin(45,25;-1,-1;5)\Lin(39,19;-1,-1;5)
\else
\Lin(57,37;-1,-1;4)\Lin(52,32;-1,-1;4)\Lin(47,27;-1,-1;4)
\Lin(37,17;-1,-1;4)\Lin(32,12;-1,-1;4)\Vec(27, 7;-1,-1;4)
\fi
\Lin(33,13;-1,-1;5)\Vec(27, 7;-1,-1;5)
\Put(26,14){\al\<{\A,\B,\C}}
%
\Lin(60,37;0,-1;3)\Lin(60,32;0,-1;3)\Lin(60,27;0,-1;3)\Lin(60,22;0,-1;3)
\Lin(60,17;0,-1;3)\Lin(60,12;0,-1;3)\Vec(60, 7;0,-1;3)
\Put(62,20){t}
}
\EndDiag%
Here $t=(\nu\<{\A,\B}\.\pi\<{\A\#\B,\C})*\nu\<{\A\#\B,\C}$

Then for any object $\D$ and for any morphisms $a\:\D\>\A$, $b\:\D\>\B$, and $c\:\D\>\C$ we
have
\begin{equation} \lf{9}
  \al\<{\A,\B,\C}\.\((a*b)*c\) = a*(b*c).
\end{equation}

\subsection{Elementary axioms for categories of information transformers} %
In this subsection we set forward the main properties of categories of ITs. All the following study
will rely exactly on these properties.

In{}~%
\ct{Gol:MSAlg, Gol:InfCatLinSys, Gol:RelInf, Gol:MV_Inf, Gol:CatInfTran, Gol:MS}
it is shown  that classes of information transformers can be
considered as morphisms in certain categories. As a rule, such categories do not have products,
which is a peculiar expression of nondeterministic nature of ITs in these categories. However, it
turns out that deterministic information transformers, which are usually determined in a natural way
in any category of ITs, form a subcategory with products. This point makes it possible to define a
``product'' of objects in a category of ITs. Moreover, it provides an axiomatic way to describe an
extension of the product operation from the subcategory of deterministic ITs to the whole category
of ITs.

\begin{definition} \lp{IT:Simple}
We shall say that a category $\CC$ is a \NW{category of information transformers} if the following
axioms hold:
\begin{enumerate}
\item  \label{SubDet}
There is a fixed subcategory of \NW{deterministic} ITs $\DD$ that contains all the objects of the
category $\CC$ ($\Ob(\DD)=\Ob(\CC)$).
\item
The classes of \NW{isomorphisms} in $\DD$ and in $\CC$ coincide,
that is, all the isomorphisms in $\CC$ are deterministic.
\item
The categories $\DD$ and $\CC$ have a common \NW{terminal object} $\Z$.
\item
The category $\DD$ has pairwise \NW{products}.
\item \label{GenProd}
There is a specified \NW{extension of morphism product} from the
subcategory $\DD$ to the whole category $\CC$, that is, for any
object $\D$ and for any pair of morphisms $a\:\D\>\A$ and
$b\:\D\>\B$ in $\CC$ there is certain information transformer
$a*b\:\D\>\A\#\B$ (which is also called a \NW{product} of ITs $a$
and $b$) such that
\[
  \pi\<{\A,\B}\.(a*b)=a, \qquad
  \nu\<{\A,\B}\.(a*b)=b.
\]
\item \label{GenParProd}
Let $a\:\A\>\C$ and $b\:\B\>\D$ are arbitrary ITs in $\CC$, then
the IT $a\*b$ defined by Eq.{}~\rf{3} satisfies Eq.{}~\rf{5}:
\[
  (a\*b)\.(c*d) = (a\.c)*(b\.d).
\]
\item \label{ProdIso}
Equality{}~\rf{7} holds not only in $\DD$ but in $\CC$ as well, that is, \NW{product} of
information transformers is  ``\NW{commutative} up to isomorphism.''
\item \label{AssocIso}
Equality{}~\rf{9} also holds in $\CC$. In other words, \NW{product} of information transformers
is ``\NW{associative} up to isomorphism'' too.
\end{enumerate}
\end{definition}

Now let us make several comments concerning the above definition.

We stress that in the description of the extension of morphism
product from the category $\DD$ to $\CC$ (Axiom~\ref{GenProd}) we
{\em do not require the uniqueness} of an IT $c\:\D\>\A\#\B$ that
satisfies conditions (1).

Nevertheless, it is easily verified, that the equations{}~\rf{4} are valid for $c=a\*b$ not only in the
category $\DD$, but in $\CC$ as well, that is,
\[
  \pi\<{\C,\D}\.(a\*b) = a\.\pi\<{\A,\B}, \qquad
  \nu\<{\C,\D}\.(a\*b) = b\.\nu\<{\A,\B}.
\]
However, the IT $c$ that satisfy the equations{}~\rf{4} may be not unique. Note also that in the
category $\CC$ Eq. (2) in general does not hold.

Further, note that Axiom{}~%
 \ref{GenParProd}
immediately implies
\[
  (a\*b)\.(c\*d) = (a\.c)\*(b\.d).
\]

Finally, note that any category that has a terminal object and
pairwise products can be considered as a category of ITs in which
all information transformers are deterministic.

\section{Fuzzy logic and fuzzy quantifiers}

 In this section we shall introduce the fuzzy interpretation of 
formulas of first order logic which is an extension of the classical 
interpretation.

\subsection{Truth values and operations of fuzzy logic}

Let $\I$ be the closed interval $\I = [\0, \1]$ -- the set of {\it 
fuzzy truth values}.  For the sake of simplicity of notations we shall 
denote any fuzzy proposition a and its truth value (evaluation of $a$) 
by the same symbol, i.e., $a\in \I$, but we shall mark fuzzy logical 
operations by dot in order to avoid ambiguity.

Let $a$ and $b$ be any fuzzy propositions, i.e., $a,b \in \I$. We put 
by definition for {\it fuzzy disjunction, conjunction, negation}, and 
{\it implication}:

\beq
a\, \dotvee \, b \, 
 &\!\!\!\ddef\!\!\!& 
\, \max (a,b), \qquad
\dot{\urcorner} \, a \ddef\, 1-a, \\
 a\dotwedge\, b \, 
 &\!\!\!\ddef\!\!\!& 
\min(a,b), \qquad a \dotar\, b\, \ddef\, \dot{\urcorner}\,a\,  \dotvee\, b.
\eeq

Let $A(x)$ be any $x$-indexed family of propositions $x\in \clX$, then 
define {\it multivalent disjunction} and {\it conjunction}:

\beq
\dotbigvee_x\, A(x)\, \ddef\, \sup_x A(x), \qquad \dotbigwedge_x\, A(x) 
\, \ddef \, \inf_x A(x).
\eeq

\subsection{Fuzzy sets as fuzzy predicates}

When we consider some family of fuzzy propositions $A(x)$, $x\in \clX$ 
we may say that $A$ is a fuzzy property (predicate, characteristics) 
for elements of the set $\clX$. That is an element $x$ of $\clX$ 
satisfies $A$ with the degree $A(x)\in \I$. On the other hand any fuzzy 
property of elements of $\clX$ may be considered as a fuzzy subset of 
$\clX$.

 In what follows we shall identify a fuzzy property, the corresponding 
fuzzy set and its characteristic function and hence use the same 
notation for all three (conceptually different!) notions. We shall also 
find it very convenient to denote the grade of membership of $x$ to a 
fuzzy set $A$ by\, $\dotin A$, that is $x\,\,\dotin A\, \ddef\, A(x)$.

\subsection{Fuzzy quantifiers}

 Quantifier in classical logic may be treated as conjunction or 
disjunction of a family of propositions. Specifically, let $A(x)$ be 
any (classical) predicate, then the formula $\forall x\, A(x)$ takes the 
truth value 1 (true) iff $A(x) = 1$ for all $x$.

 Now for any fuzzy predicate $A(x)$, $x\in\clX$ define {\it fuzzy 
universal} and {\it existential quantifiers}:
\beq
\dotal\ x\,\, A(x)\, \ddef\, \dotbigwedge_x A(x), \qquad \dotex \ x \, 
A(x) \, \ddef\, \dotbigvee_x A(x).
\eeq

\subsection{Bounded quantifiers}

 Bounded quantifiers are convenient for manipulations with the 
statements of the kind: ``all $x$ such that... satisfy...'' and ``there 
exists $x$ such that... satisfying...''. We shall need fuzzy quantifiers 
with {\it fuzzy bounds}.

By analogy with the classical logic let us put by definition:
\beq
\dotal\  x\,\, \dotin\, A\,\, B(x)  &\ddef&  \dotal \ x \ 
\left( x \ \dotin \ A \dotar \ B(x)\right),\\ 
\dotex\ x \ \dotin \ A \,\, B(x)  
&\ddef&  \dotex \ x \ \left( x \ \dotin \ A\, \dotwedge \, 
B(x)\right).  
\eeq

Note that the bounded quantifiers are related in the same way as the 
unbounded ones:
\beq
\dot{\urcorner}  \left(\,  \dotal \ x \ \dotin \ A\, B(x) \right) 
=\dotex\ x\,\,\,  \dotin \ A\,\, \dot{\urcorner} B(x).  \eeq

\section{Algebra of fuzzy sets}

 The algebra of sets in the classical set theory is determined by the 
underlying classical logic.  The systematic extension of set theoretic 
notations allows to introduce the algebra of fuzzy sets in a very 
natural way.

\subsection{Definitions of fuzzy sets by comprehension}

Let $\varphi$ be any fuzzy property for elements $x$ of $\clX$. This 
property may be any formula in the fuzzy logic discussed above. By 
$A=\{x|\varphi(x)\}$ we shall denote the fuzzy set such that the level 
of membership of $x$ to $A$ is equal to the truth value of 
$\varphi(x)$, i.e.:
\be
A= \{x|\varphi(x)\} \qquad   {\rm iff} \qquad  x \ \dotin\, 
A=\varphi(x).  \ee  

\subsection{Fuzzy sets operations}

Suppose that $A$ and $B$ are any fuzzy sets over the space $\clX$. The 
principle (10) provides definitions of {\it union, intersection}, and 
{\it complement} for fuzzy sets:
\beq
A \ \dotcup \ B  &\ddef&  \{ x| x \ \dotin \ A\,\, \dotvee \ 
x \ \dotin \ B\},\\
A\dotcap \ B &\ddef& \{x| x \ \dotin \ A\dotwedge \
x \ \dotin \ B\},\\
\overline{A}  &\!\!\!\ddef\!\!\!&  \{x| \ \dot{\urcorner}\, (x \ \dotin 
\ A\}.  \eeq

Now let $\alpha(y)$ be any family of fuzzy sets on $\clX$ indexed by 
$y\in\clY$.  By definition, put
\beq
\dotbigcup_y \alpha(y)  &\ddef&  \{x| \ \dotex \  y\, x \ \dotin \ 
\alpha(y)\},\\
\dotbigcap_y \alpha(y)  &\ddef&  \{x| \ \dotal \ y\, x \ \dotin 
\alpha(y)\}.
\eeq

Finally, if the index $y$ itself varies in a fuzzy set $Y$, then 
operations for fuzzy families of fuzzy sets are easily described by 
bounded fuzzy quantifiers:
\beq
\dotbigcup_{\sss y \ \dotin Y}\, \alpha(y)  &\!\!\!\ddef\!\!\!&  
\{x| \ \dotex  \ y \ \dotin \ Y \, x \ \dotin \, \alpha(y)\},\\ 
\dotbigcap_{\sss y\ \dotin Y}\, \alpha(y)  &\!\!\!\ddef\!\!\!&  
\{x|\ \dotal \ y \ \dotin Y\, x \ \dotin \ \alpha(y)\}.  \eeq

\subsection{Containment of fuzzy sets}

In contrast to the traditional definition of containment of fuzzy sets: 
``$A\subset B$ iff $A(x)B(x)$ for all $x$'' we shall adopt another 
definition of {\it containment}, which is the natural extension of 
containment of crisp sets. In our approach ``$A$ is contained in $B$'' is 
a {\it fuzzy property}. Hence it may take fuzzy truth values. By 
definition, put
\beq
A \ \dotsubset \ B \ \ddef \ \dotal \ x \ \dotin \ A\,\, x \ \dotin \, B,
\eeq
that is ``$A$ is {\it contained} in $B$'' with the degree in which 
``{\it all} the  elements {\it contained} in $A$ {\it belong} to $B$'' 
(in these phrases we emphasized fuzzy notions).

\section{Fuzzy sets and distributions}

We shall say that a {\it fuzzy distribution} (or {\it possibility 
distribution} \ct{Dub:FSII}) $X$ on the space $\clX$ is any fuzzy set in $\clX$. 
Note that we consider a fuzzy distribution as an analog of a 
probabilistic distribution.

\subsection{Joint distributions}

We shall say that any fuzzy distribution $C$ on $\clX \times \clY$ 
determines a {\it fuzzy joint distribution} \ct{Manes:ClassFuzzy} of the pair $\langle 
x,y\rangle$, $x\in\clX$, $y\in\clY$.

 Sometimes it is convenient to interpret a joint distribution as a 
fuzzy relation, i.e., to consider $C(x,y)$ as the ``degree in which $x$ 
and $y$ satisfy the relation $C$''.

The distribution $C^\bullet$ on $\clY\times\clX$ is called the 
{\it converse} of $C$ if
\beq
C^\bullet(y,x) =C(x,y).
\eeq

The joint distribution $C$ induces the {\it marginal distributions} \ct{Gol:FuzDec, Gol:MSI}  
$X$ and $Y$ on the spaces $\clX$ and $\clY$ respectively:
\beq
X=\{x| \ \dotex \  y\, C(x,y)\}, \qquad Y=\{y|\ \dotex \  x \, C(x,y)\}.
\eeq

\subsection{Transition distributions}

 The notion of a transition distribution \ct{Barra} is of great importance in 
mathematical statistics.

We shall say that a {\it fuzzy transition distribution} from $\clX$ to 
$\clY$ is any map $\alpha$ from $\clX$ to the family of all fuzzy 
distributions on $\clY$, that is, $\alpha$ takes each element 
$x\in\clX$ to some fuzzy distribution $\alpha(x)$ on the space $\clY$.

Note that the condition $A(x,y) = \alpha(x)(y)$ for all $x$ in $\clX$ 
for all $y$ in $\clY$ determines the one-to-one correspondence between 
$A$ and $\alpha$.  Using set theoretic notations we may write:
\beq
A=\{\langle x,y\rangle| \alpha(x)(y)\}, \qquad
\alpha(x)=\{y| A(x,y)\}.
\eeq

By $\overline{\alpha}$ we shall denote the transition distribution, 
such that $\overline{\alpha}(x)\, \ddef\, \overline{\alpha(x)}$. 
Obviously, $\overline{\alpha}$ corresponds to $\overline{A}$.

 The correspondence between joint distributions and transition 
distributions reveals the {\it algebraic semantics} of the marginal 
distributions. Let $\alpha$ and $\beta$ be transition distributions 
corresponding to $A$ and $A^\bullet$ respectively. Then
\beq
X  &\!\!\!=\!\!\!&  \{x| \ \dotex \ y \,\,\, \beta(y)(x)\} =\dotbigcup_y \
\beta(y),\\
Y &\!\!\!=\!\!\!& \{y| \ \dotex \ x \,\,\, \alpha(x)(y)\} =\dotbigcup_x \
\alpha(x).
\eeq

\subsection{Images of distributions}

 As a probabilistic transition distribution transforms one probability 
distribution to another \ct{Barra}, in the same way a fuzzy transition 
distribution takes some fuzzy distribution to another one.

Let $X$ be any distribution on $\clX$ and a is some transition 
distribution from $\clX$ to $\clY$. We say that a distribution 
$\alpha\clX$ on $\clY$ is the {\it image of the distribution $X$ 
induced by the transition distribution} $\alpha$ if
\beq
\alpha X\, \ddef \,  \{y| \ \dotex \, x \,\, \dotin  X\, y \ \dotin \
\alpha(x)\} = \dotbigcup_{\sss x \ \dotin  X} \alpha(x).  \eeq

 There is an interesting notion dual to the notion of the image. We say 
that a distribution $\alpha \circ \clX$ on $\clY$ is the {\it lower 
image of the distribution $X$ induced by the transition distribution} 
$\alpha$ if
\beq
\alpha\circ X\, \ddef\, \{y| \ \dotal  x \ \dotin  X\,\,\, y \ \dotin \
\alpha(x)\} =\dotbigcap_{\sss x\ \dotin  X} \alpha(x).  \eeq

It may be proved that $\overline{\alpha\circ X} =\overline{\alpha} X$.

\subsection{Generated joint distributions}

Suppose that $X$ is an arbitrary distribution in $\clX$ and $\alpha$ is 
some transition distribution from $\clX$ to $\clY$. We say that the 
distribution $\alpha* X$ in $\clX\times \clY$ is the {\it joint 
distribution generated} by $\alpha$ and $X$ if
\beq
\alpha * X\, \ddef\, \{\langle x,y\rangle | x \ \dotin \, 
X\dotwedge \ y \,\, 
\dotin \ \alpha(x)\}.  \eeq 

It can be easily verified that the image $\alpha X$ coincides with the 
$\clY$-marginal distribution of the generated joint distribution 
$\alpha* X$.

\subsection{Fuzzy conditional distributions}

 We shall define this notion in the same way as its probabilistic 
prototype is defined in statistics  \ct{Gol:FuzDec, Gol:MSI, Barra}.

Let $A$ be any joint distribution on $\clX\times \clY$ . We say that 
some transition distribution $\alpha$ from $\clX$ to $\clY$ is (a 
variant of) the {\it conditional distribution} for $A$ with respect to 
$\clX$ if $A$ is generated by its $\clX$-marginal distribution $X$ 
together with $\alpha$:
\beq
A=\alpha * X.
\eeq

Similarly we say that $\beta$ from $\clX$ to $\clY$ is a conditional 
distribution for $A$ with respect to $\clY$ if
\beq
A^\bullet =\beta * Y.
\eeq

\begin{theorem}
For any joint distribution $A$ in $\clX \times\clY$, conditional 
distributions always exist.  Some variants of conditional distributions 
$\alpha$ and $\beta$ are determined by the following expressions:
\beq
\alpha(x) =\{x|\, A(x,y)\}, \qquad \beta(y)=\{y|\, A(x,y)\}.
\eeq
\end{theorem}

\subsection{Iterated quantifiers}

 The possibility of representation of any joint distribution by 
conditional distributions has interesting consequences that will be 
used in the Bayesian decision analysis.

\begin{theorem}
Let $A$ be any joint distribution on $\clX\times\clY$, $X$ and $Y$ 
--- its marginal distributions, $\alpha$ and $\beta$ --- conditional 
distributions, that is $A = \alpha * X  =(\beta * Y)^\bullet$. Then 
for any formula $\varphi$
\beq
\dotal \ \langle x,y\rangle \ \dotin \ A\varphi  
&\!\!\!=\!\!\!& 
\dotal  \ x \,\, \dotin \ X \ \dotal \ y \ \dotin \alpha(x) \varphi = \dotal 
\ y \ \dotin \ Y \ \dotal \ x \ \dotin \ \beta(y) \varphi,\\ 
\dotex \ \langle x,y\rangle \ \dotin \ A \varphi  
&\!\!\!=\!\!\!& 
\dotex \ x \ \dotin \ X \ \dotex \ y \ \dotin \ \alpha(x)\varphi = 
\dotex \ y \ \dotin \ Y \ \dotex \ x \ \dotin \, \beta(y)\varphi.  \eeq 
\end{theorem}

\section{Semantics of decision making problems}

\subsection{Fuzzy games with a priory information}

 Let $\clX$ be some space of objects of interest; let $\clD$ be some 
{\it space of decisions}; and let $G$ be a fuzzy relation between 
$\clX$ and $\clD$, i.e., a fuzzy set in $\clX\times\clD$ that 
determines the notion "good decision," i.e., $G =\{\langle 
x,d\rangle|d$ is a good decision for $x\}$.

 We are dealing in fact with the two-person game \ct{Black:Gir} 
$\langle\clX, \clD, G\rangle$, where $\clX$ is the decision space of the first 
player, $\clD$ is the decision space of the second player, and the 
fuzzy notion of the quality of decisions for the second player $G$ 
replaces the classical loss function.  Hence we shall call 
$\langle\clX, \clD, G\rangle$ 
the fuzzy game \ct{Gol:FuzDec, Gol:MSI}. If it is known a priori that all choices of the 
first player are restricted by a fuzzy set $X$, then we are dealing 
with the {\it fuzzy game with the a priory information} \ct{Gol:FuzDec, Gol:MSI}.

\subsection{Good decisions}

Let us now define the notion $G_X$ = ``{\it for all $x$ in} $X$ decision 
$d$ is {\it good} for $x$'' i.e., formally:
\be
G_X \, \ddef\, \dotal  \ x \ \dotin \ X\, G(x,d).
\ee
Truth value of this formula determines the quality of the decision $d$ 
for $X$. Note that the distribution of good decisions permits the 
following interpretation. Let the transition distribution $\gamma$ from 
$\clX$ to $\clD$ is determined by $\gamma(x)=\{d| G(x,d)\}$, i.e., the 
distribution of ``good decisions for $x$''.  Let $\gamma^\bullet$ be the 
opposite to $\gamma$.  Then
\beq
G_X =\{d| \ \dotal \ x \ \dotin \ X\, x \ \dotin  \ \gamma^\bullet(d)\} 
=\{x|X\, \dotsubset \ \gamma^\bullet (d)\}.  \eeq

Another interpretation for $G_X$ --- its algebraic semantics --- comes 
from the representation:
\beq
G_X= \dotbigcap_{\sss x \ \dotin X} \gamma(x) =\gamma\circ X.
\eeq
 Note that in statistics the loss function is often used. Suppose that 
some joint fuzzy distribution $B$ on $\clX\times\clD$ is fixed and let 
$\beta$ be the associated transition distribution $\beta(x)=$ ``{\it 
bad} decisions for $x$'' $= \{d|\langle x,d\rangle \,\, \dotin \ B\}$. Then 
the fuzzy set of bad decisions with respect to the a priory 
distribution $X$ may be defined by:  
\beq 
B_X =\{d| \ \dotex \ x \ \dotin \ 
X\, B(x,d)\} =\dotbigcup_{\sss x\ \dotin \ X} \beta(x) =\beta X.  \eeq

However if $B = \overline{G}$ (it is quite natural) then we have 
$\beta=\overline{\gamma}$ and
\beq
B_X =\beta X =\overline{\gamma} X = \overline{\gamma\circ X} 
=\overline{G_X},
\eeq
so the two approaches are dual to each other.

Let us consider an {\it optimal decision problem} for the second player 
with respect to the a priory information $X$. The statement ``there 
exists an optimal decision for the second player'' may be expressed by
\be
\dotex \ d \ \dotal  \ x \ \dotin \ X\, G(x,d)\, = \,\dotex  \ d \ \dotin \ 
G_X.  
\ee 
We shall say that a decision $d_X$, is {\it optimal} (or a 
{\it Bayes decision}) {\it with respect to the a priori distribution} 
$X$ (or simply $X-optimal$) if $G_X(d)$ takes on its maximum at $d_X$.

\subsection{Optimal decision strategies}

 Now consider problems of constructing optimal decision strategies in 
fuzzy experiments.  Suppose that a fuzzy experiment from $\clX$ to 
$\clY$ is determined by a fuzzy transition distribution $\alpha$. A 
{\it decision strategy} (or simply {\it strategy}) $r$ is any mapping 
from the observation space $\clY$ to the decision space $\clD$.

Let us define the relation $H(x, r) =$ ``decision strategy $r$ is good 
for $x$'' by
\beq
H(x,r)\, \ddef\, \dotal \ y \ \dotin \ \alpha(x)\, r(y) \,\, \dotin \
\gamma(x).  \eeq It may be shown that $H(x,r) = r\alpha(x)\,\dotsubset 
\gamma(x)$.

In fact, we have come to a new two-person game 
$\langle\clX, \clD^\clY,H \rangle$, where the decision space
 of the second player is replaced by the set of all strategies, and the 
goodness relation $G$ for decisions is replaced by the goodness for 
strategies. Such games are analogous to classical statistical games 
\ct{Black:Gir}.

Finally for a given a priory distribution $X$ we define the {\it Bayes 
goodness} for a strategy with respect to $X$, i.e., the fuzzy set on 
$\clD^\clY$, $H_X(r)=$ \,  ``$r$ is good for all $x$ in $X$''
\be
H_X(r)\, \ddef\, \dotal \ x \ \dotin \ X\, H(x,r).
\ee

The {\it optimal decision strategy problem} is one of obtaining a 
mapping $r_X$, called a {\it Bayes strategy}, such that the level of 
goodness for $r_X$ in (13) is the highest.

\subsection{Bayes principle}

 A well-known result of the theory of statistical games states that an 
optimal decision strategy problem for a statistical game may be reduced 
to a similar problem for an original game \ct{Black:Gir}.  This reduction involves 
conditional distributions.

\begin{theorem}
Let $X$ be any a priori distribution in $\clX$, $\alpha$ any transition 
distribution from $\clX$ to $\clY$, and $\beta$ a variant of a 
conditional distribution for $\alpha * X$ with respect to $\clY$. 
Assume that for every $y\, in\, $$\clY$  there exists a Bayes decision 
$d_{\beta(y)}$, which is optimal for the original decision problem with 
respect to the distribution $\beta(y)$.  Then the optimal decision 
strategy $r_X$ and the corresponding goodness $H_X(r_X)$ are determined 
by
\beq
r_X(y) =d_{\beta(y)}, \qquad H_X(r_X) =\dotal \ y \ \dotin \ Y\, 
G_{\beta(y)} \left(d_{\beta(y)}\right).
\eeq

\end{theorem}

For an arbitrary decision strategy $r$
\beq
H_X(r)  &\!\!\!=\!\!\!& 
\dotal  \ x \ \dotin \ X \ \dotal \ y \ \dotin \ \alpha X\, r(y)\,\, \dotin \
\gamma(x) = \dotal \ \langle x,y\rangle  \ \dotin \ \alpha * X\, r(y)\, 
\, \dotin \ \gamma(x)=\\ 
&\!\!\!=\!\!\!& 
\dotal  \ y \ \dotin \ \alpha 
X\,\, \dotal \ x \ \dotin\  \beta Y\, r(y)\,\, \dotin \ \gamma(x) = \dotal 
\ y \ \dotin \ \alpha X\, r(y)\,\, \dotin\, G{\beta(y)}.
 \eeq
 In the calculations above we utilize the definition of the conditional 
distribution and in the last step we use the definition of $G_X(d)$ (2) 
for $X = \beta(y)$.

Hence, by virtue of the optimality of the Bayes decision 
$d_{\beta(y)}$, we have for any $y$
\beq
G_{\beta(y)}\left( r(y)\right) G_{\beta(y)} \left( d_{\beta(y)}\right) 
=G_{\beta(y)} \left( r_X(y)\right).
\eeq
Therefore we obtain
\beq
H_X(r) = \dotbigvee_{\sss y\ \dotin  \alpha X}\, G_{\beta(y)}\, \left( 
r(y)\right) \dotbigvee_{\sss y \ \dotin  \alpha X}\, G_{\beta(y)}\, 
\left( r_X(y)\right) =H_X(r_X).  \eeq

 Note that the notion of conditional distribution and tightly related 
with it ``interchange'' of quantifiers play the principal role in the 
proof.

 The Bayes principle is very intuitive. It asserts that to construct 
(or calculate) an optimal strategy $r_X$ for a given observation $y$ 
one has to find the conditional distribution of $x$ for a given $y$, 
i.e., $\beta y$ and then take a decision $d_{\beta(y)}$, which is 
optimal with respect to this distribution. In other words, the 
observation of $y$ results in the passage from the a priori information 
$X$ to the a posteriory information $\beta(y)$.

\section{Informativeness of information transformers}

\subsection{Accuracy relation} 
In order to define informativeness relation we will need to
introduce first the following auxiliary notion.

\begin{definition}
We will say that $\HP$ is an \NW{accuracy} relation on an
IT-category $\CC$ if for any pair of objects $\A$ and $\B$ in
$\CC$ the set $\CC(\A,\B)$ of all ITs from $\A$ to $\B$ is
equipped with a partial order $\HP$ that satisfies the following
\NW{monotonicity} conditions:
\[
  a \HP a',\; b \HP b' \Imp a\.b \HP a'\.b',
\]
\[
  a \HP a',\; b \HP b' \Imp a*b \HP a'*b'.
\]
\end{definition}

Thus, the composition and the product are monotonous with respect
to the partial order $\HP$. For a pair of ITs $a,b\in\CC(\A,\B)$
we shall say that $a$ is more \NW{accurate} then $b$ whenever
$a\HP b$.

It obviously follows from the very definition of the operation
$\*$~\rf{3} and from the monotonicity conditions that the
operation $\*$ is monotone as well:
\[
  a \HP a',\; b \HP b' \Imp a\*b \HP a'\*b'.
\]

It is clear that for any IT-category there exists at least a
``trivial variant'' of the partial order $\HP$, namely, one can
choose an equality relation for $\HP$, that is, one can put $a\HP
b \iffdef a=b$. However, many categories of ITs (for example,
multivalued and fuzzy ITs) provide a ``natural'' choice of the
accuracy relation, which is different from the equality relation.

\subsection{Definition of informativeness relation} 
Suppose $a\:\D\>\A$ and $b\:\D\>\B$ are two information transformers with a common source
$\D$. Assume that there exists an IT $c\:\A\>\B$ such that $c\.a=b$. Then any information that
can be obtained from $b$ can be obtained from $a$ as well (by attaching the IT $c$ next to $a$).
Thus, it is natural to consider the information transformer $a$ as being more informative than the IT
$b$ and also more informative than any IT less accurate than $b$.

Now we give the formal definition of the informativeness relation in the category of information
transformers.

\begin{definition}  
We shall say that an information transformer $a$ is more \NW{informative} (better) than $b$ if there
exists an information transformer $c$ such that $c\.a\HP b$, that is,
\[
  a \BI b \Iffdef \exists c \q c\.a\HP b.
\]
\end{definition}

It is easily verified that the informativeness relation $\BI$ is a preorder on the class of information
transformers in $\CC$. This preorder $\BI$ induces an equivalence relation $\EI$ in the following
way:
\[
  a \EI b \Iffdef a \BI b \AND b \BI a.
\]

Obviously, the relation ``more informative'' extends the relation ``more accurate,'' that is,
\[
  a \HP b \imp a \BI b.
\]

\subsection{Main properties of informativeness}
It can be easily verified that the informativeness relation $\BI$ satisfies the following natural
properties.

\begin{lemma} 
Consider all information transformers with a fixed source $\D$.

\begin{itemize}
\item[\It{a}] The identity information transformer $i\<\D$ is the most informative and the
terminal information transformer $z\<\D$ is the least informative:
\[
  \forall a \quad i\<\D \BI a \BI z\<\D.
\]
\item[\It{b}] Any information transformer $a\:\D\>\spc{B}\#\spc{C}$ is more informative
than its parts $\pi\<{\B,\C}\.a$ and $\nu\<{\B,\C}\.a$.
\item[\It{c}] The product $a*b$ is more informative than its components
\[
  a*b \BI a,b.
\]
\end{itemize}
\end{lemma}

Furthermore, the informativeness relation is compatible with the composition and the product
operations.

\begin{lemma} 
\It{a} If $a \BI b$, then $a\.c \BI b\.c$.

\It{b} If $a \BI b$ and $c \BI e$, then $a*c \BI b*e$.
\end{lemma}

\subsection{Structure of the family of informativeness equivalence classes}
Let $a$ be some information transformer. We shall denote by $[a]$ the equivalence (with respect to
informativeness) class of $a$. We shall also use boldface for equivalence classes, that is, $a\in\EC
a$ is equivalent to $\EC a=[a]$.

\begin{theorem} 
Let $\mathfrak{J}(\D)$ be the family of informativeness equivalence classes for the class of all
information transformers with a fixed domain $\D$. The family $\mathfrak{J}(\D)$ forms a
partial ordered Abelian monoid $\TPL{\mathfrak{J}(\D),\BI,*,\EC 0}$ with the smallest
element $\EC 0$ and the largest element $\EC 1$, where
\LongEq{
  [a] \BI [b] \iffdef a \BI b,
\qquad
  [a] * [b] \Eqdef [a*b],
}{
  \EC 0 \Eqdef [z\<\D],
\qquad
  \EC 1 \Eqdef [i\<\D].
}

Moreover, the following properties hold:

\begin{itemize}
\item[\It{a}]  $\quad \EC {0*a=a}$,
\item[\It{b}]  $\quad \EC {1*a=1}$,
\item[\It{c}]  $\quad \EC {0 \LI a \LI 1}$,
\item[\It{d}]  $\quad \EC {a*b \BI a,b}$,
\item[\It{e}]  $\quad \EC {(a \BI b) \AND (c \BI e) \Imp a*c \BI b*e}$.
\end{itemize}
\end{theorem}

\subsection{Informativeness structure of concrete IT-categories}
\subsubsection{Linear stochastic ITs}
\begingroup
\DefLin
In addition to the trivial relation of accuracy (which coincides with the equality relation) one can
define the accuracy relation in the following way:
\[
  \TPL{A_a,\S_a}\HP\TPL{A_b,\S_b} \Iffdef A_a=A_b,\;\S_a\LE\S_b.
\]
However, it can be proved that the informativeness relations corresponding these different accuracy
preorders, actually coincide.

\begin{theorem}
In the category of linear information transformers every equivalence class $[a]$ corresponds to a
pair $\TPL{\spc{Q},S}$,
where $\spc{Q}\SS\D$ is an Euclidean subspace and $S\:\spc{Q}\>\spc{Q}$ is nonnegative
definite operator, that is, $S\GE0$. In these terms
\[
  \TPL{\spc{Q}_1,S_1}\GE\TPL{\spc{Q}_2,S_2}
  \Iffdef
  \spc{Q}_1\supseteq\spc{Q}_2, \q S_1\upharpoonright\spc{Q}_2\LE S_2.
\]
Here $S_1\upharpoonright\spc{Q}_2$ (the restriction of $S_1$ on
$\spc{Q}_2$) is defined by the expression
$S_1\upharpoonright\spc{Q}_2 \Eqdef P_2I_1S_1P_1I_2$, where
$I_j\:\spc{Q}_j\>\D$ is the subspace inclusion, and
$P_j\:\D\>\spc{Q}_j$ is the orthogonal
projection (cf.{}~%
\ct{Gol:InfCatLinSys, Pyt:ZadRed}).
\end{theorem}
\endgroup

\subsubsection{The category of sets as a category of ITs}
It is not hard to prove that for a given set $\D$, the class of equivalent informativeness for an IT $a$
with the set $\D$ being its domain, is completely determined by the following equivalence relation
$\EQ_a$ on $\D$:
\[
  x \EQ_a y \iffdef ax = ay \qquad \forall x,y\in\D.
\]
Furthermore, $a\BI b$ if and only if the equivalence relation $\EQ_a$ is \NW{finer} than $\EQ_b$,
that is,
\[
  a\BI b \Iff \forall x,y\in\D \;\; \(x \EQ_a y \imp x \EQ_b y\).
\]

It is clear that we have the following

\begin{theorem}
The partially ordered monoid of equivalence classes for ITs with the source $\D$, is isomorphic
to the monoid of all equivalence relations on $\D$ equipped with the order ``finer'' and with the
product:
\[
  x\; (\EQ_a*\EQ_b)\; y \Iffdef \(x \EQ_a y,\; x \EQ_b y\) \qquad
  \forall x,y\in\D.
\]
\end{theorem}

\subsubsection{Multivalued ITs}
\begingroup
\DefMul
In addition to the trivial accuracy relation in the category of multivalued ITs one can put
\[
  a \HP b \Iffdef \forall x\in\D\; ax\SS bx.
\]
These two accuracy relations lead to different informativeness relations{}~%
\ct{Gol:MV_Inf},
called (strong) informativeness $\BI$ and weak informativeness $\BIW$.

For the both informativeness relations the classes of equivalent ITs with a fixed source $\D$ can be
described explicitly.

\begin{theorem}
In the category of multivalued ITs with weak informativeness
every class of equivalent ITs corresponds to a certain covering
$\P$ of the set $\D$, such that if $\P$ contains some set $B$ then
it contains all its subsets:
\[
    \(\exists B\in\P\;(A\SS B)\) \Imp A\in\P.
\]
Moreover, a covering $\P_1$ is more (weakly) informative than $\P_2$ (namely, $\P_1$
corresponds to a class of more (weakly) informative ITs than $\P_2$) if $\P_1$ is contained in
$\P_2$, that is,
\[
  \P_1 \BIW \P_2 \Iffdef \P_1 \SS \P_2.
\]
\end{theorem}

\begin{theorem}
In the category of multivalued ITs with strong informativeness
every class of equivalent ITs corresponds to a covering $\P$ of
the set $\D$, that satisfy the following condition:
\[
  \(\Big.
    \(\exists B\in\P\;A\SS B\)
    \;\&\;
    \(\exists \spc{B}\SS\P\;A=\cUp\spc{B}\)
  \)
\iftc
  \imp A\in\P.
\else
  \Imp A\in\P.
\fi
\]
In this case
\begin{eqnarray*}
  \P_1 \BI \P_2
  &\iffdef&
  \Big(
    \(\forall A\in\P_1 \q \exists B\in\P_2 \q A \SS B\)
  \\ &&
    \;\&\;
    \(\forall B\in\P_2 \q \exists \spc{A}\SS\P_1 \q B=\cUp\spc{A}\)
  \Big).
\end{eqnarray*}
\end{theorem}
\endgroup

\section{Informativeness and decision-making problems}
\ls{SemInf}
In this section, we consider an alternative (with respect to the above) approach to informativeness
comparison. This approach is based on treating information transformers as data sources for
decision-making problems.

\subsection{Decision-making problems in categories of ITs}
Results of observations, obtained on real sources of information (e.g. indirect  measurements) are as
a rule unsuitable for straightforward interpretation. Typically it is assumed that observations suitable
for interpretation are those into a certain object $\U$ which in what follows will be called object of
\NW{interpretations} or object of \NW{decisions}.

By an \NW{interpretable} information transformer for signals from an object $\D$ we mean any
information transformer $a\:\D\>\U$.

It is usually thought that some interpretable information transformers are more suitable for
interpretation (of obtained results) than others. Namely, on a set $\CC(\D,\U)$ of information
transformers from $\D$ to $\U$, one defines some preorder relation $\BQ$, which specifies the
relative quality of various interpretable information transformers. Typically the relation $\BQ$ is
predetermined by the specific formulation of a problem of optimal information transformer synthesis
(that is, decision-making problem).

We shall say that an abstract \NW{decision-making problem} is determined by a triple
$\TPL{\D,\U,\BQ}$, where $\D$ is an object of studied (input) signals, $\U$ is an object of
decisions (or interpretations), and $\BQ$ is a preorder on the set $\CC(\D,\U)$.

We shall call a preorder $\BQ$ \NW{monotone} if for any $a,b\in\CC(\D,\U)$
\[
  a\HP b \imp a\BQ b,
\]
that is, more accurate IT provides better quality of interpretation.

For a given information transformer $a\:\D\>\A$ we shall also say that an IT $b$ \NW{reduces}
$a$ to an interpretable information transformer if $b\.a\:\D\>\U$, that is, if $b\:\A\>\U$. Such
an information transformer $b$ will be called a \NW{decision strategy}.

The set of all interpretable information transformers obtainable on the basis of $a\:\D\>\A$ will be
denoted $\UU_a\SS\CC(\D,\U)$:
\[
  \UU_a \Eqdef \SET{b\.a|b\:\A\>\U}.
\]

We shall call a decision strategy $r\:\A\>\U$ \NW{optimal} (for the IT $a$ with respect to the
problem $\TPL{\D,\U,\BQ}$) if the IT $r\.a$ is a maximal element in $\UU_a$ with respect to
$\BQ$. Thus, a decision-making problem for a given information transformer $a$ is stated as the
problem of constructing optimal decision strategies.

\subsection{Semantical informativeness}
The relation $\BQ$ induces a preorder relation $\BS$ on a class of information transformers
operating from $\D$ in the following way.

Assume that $a$ and $b$ are information transformers with the source $\D$, that is, $a\:\D\>\A$,
$b\:\D\>\B$. By definition, put
\[
  a \BS b \Iffdef
  \forall b'\:\B\>\U \q \exists a'\:\A\>\U \q
  a'\.a \BQ b'\.b.
\]
In other words, $a \BS b$ if for every interpretable information transformer $d$ derived from $b$
there exists an interpretable information transformer $c$ derived from $a$ such that $c \BQ d$, that
is,
\[
  a \BS b \Iff \forall d\in\UU_b \q \exists c\in\UU_a \q c \BQ d.
\]

It can easily be checked that the relation $\BS$ is a preorder relation.

It is natural to expect that if one information transformer is more informative than the other, then the
former will be better than the latter in any context. In other words, for any preorder $\BQ$ on the set
of interpretable information transformers the induced preorder $\BS$ is dominated by the
informativeness relation $\BI$ (that is, $\BS$ is weaker than $\BI$). The converse is also true.

\begin{definition} 
We shall say that an information transformer $a$ is \NW{semantically more informative} than $b$ if
for any interpretation object $\U$ and for any preorder $\BQ$ (on the set of interpretable
information transformers) $a\BS b$ for the induced preorder $\BS$.
\end{definition}

The following theorem is in some sense a ``completeness'' theorem, which establishes a relation
between ``structure'' ($b$ can be ``derived'' from $a$) and ``semantics'' ($a$ is uniformly better then
$b$ in decision-making problems).

\begin{theorem} 
For any information transformers $a$ and $b$ with a common source $\D$, information transformer
$a$ is more informative than $b$ if and only if $a$ is semantically more informative than $b$.
\end{theorem}

Let us remark that the above proof relies heavily on the extreme extent of the class of decision
problems involved. This makes it possible to select for any given pair of ITs $a,b$ an appropriate
decision-making problem $\TPL{\D,\U_b,\BQ_b}$ in which the interpretation object $\U_b$ and
the preorder $\BQ_b$ depend on the IT $b$. However, in some cases it is possible to point out a
concrete (universal) decision-making problem such that
\[
  a\BI b \Iff a\BS b.
\]

\begin{theorem} 
Assume that for a given object $\D$ there exists an object $\widetilde\D$ such that for every
information transformer acting from $\D$ there exists an equivalent (with respect to informativeness)
IT acting from $\D$ to $\widetilde\D$, that is,
\[
  \forall\B \q \forall b\:\D\>\B \q \exists b'\:\D\>\widetilde\D \q b\EI b'.
\]
Let us choose the decision object $\U \Eqdef \widetilde\D$ and the preorder $\BQ$, defined by
\[
  c\BQ d \Iffdef c\HP d.
\]
Then $a\BI b$ if and only if $a\BS b$.
\end{theorem}

Note that in general case an optimal decision strategy (if exists) can be nondeterministic. However,
in many cases it is sufficient to search optimal strategies among deterministic ITs. Indeed, in some
categories of information transformers the relation of ``accuracy'' satisfies the following condition:
every IT is dominated by some deterministic IT, that is, for every IT there exists a more accurate
deterministic IT.

\begin{proposition} \lp{DetMajor}
Assume that $\TPL{\D,\U,\BQ}$ is a monotone decision-making problem in a category of ITs
$\CC$. Assume also that the following condition holds:
\[
  \forall c\in\Ar(\CC)\q \exists d\in\Ar(\DD)\q d \HP c.
\]
Then for any IT $a\:\D\>\R$ and for any decision strategy $r\:\R\>\U$ there exists a
deterministic strategy $r\<0\:\R\>\U$ such that $r\<0\.a \BQ r\.a$.
\end{proposition}

\subsection{Decision-making problems with a prior information}
In this section we formulate in terms of categories of information transformers an analogy for the
classical problem of optimal decision strategy construction for decision \NW{problems with a prior
information} (or information \NW{a priori}). We also prove a counterpart of the \NW{Bayesian
principle} from the theory of statistical games{}~%
\ct{Borovkov, Black:Gir}.
Like its statistical prototype it reduces the problem of constructing an optimal decision strategy to a
much simpler problem of finding an optimal decision for a \NW{posterior information} (or
information \NW{a posteriori}).

First we define in terms of categories of information transformers some necessary concepts, namely,
concepts of distribution, conditional information transformer, decision problem with a prior
information, and others.

\subsection{Distributions in categories if ITs}
\ls{Distr}%
We shall say that a \NW{distribution} on an object $\A$ (in some fixed category of ITs $\CC$) is
any IT $f\:\Z\>\A$, where $\Z$ is the terminal object in $\CC$.

The concept of distribution corresponds to the general concept of an element of some object in a
category, namely, a morphism from the terminal object (see, e.g.,{}~%
\ct{Goldblatt}).

Any distribution of the form $h\:\Z\>\A\#\B$ will be called a \NW{joint distribution} on $\A$ and
$\B$. The projections $\pi\<{\A,\B}$ and $\nu\<{\A,\B}$ on the components $\A$ and $\B$
respectively, ``extract'' \NW{marginal distributions} $f$ and $g$ of the joint distribution $h$, that is,
\[
  f=\pi\<{\A,\B}\.h\:\Z\>\A,
\]
\[
  g=\nu\<{\A,\B}\.h\:\Z\>\B.
\]

We say that the components of a joint distribution $h\:\Z\>\A\#\B$ are \NW{independent}
whenever this joint distribution is completely determined by its marginal distributions, that is,
\[
  h=(\pi\<{\A,\B}\.h)*(\nu\<{\A,\B}\.h).
\]

Let $f$ be an arbitrary distribution on $\A$ and let $a\:\A\>\B$ be some information transformer.
Then the distribution $g=a\.f$ in some sense ``contains an information about $f$.'' This concept
can be expressed precisely of one consider the joint distribution \NW{generated} by the distribution
$f$ and the IT $a$:
\[
  h\:\Z\>\A\#\B, \qquad h=(i\<\A*a)\.f.
\]
Note, that the marginal distributions for $h$ coincide with $f$ and $g$, respectively. Indeed,
\[
  \pi\<{\A,\B}\.h = \pi\<{\A,\B}\.(i\<\A*a)\.f = i\<\A\.f = f,
\]
\[
  \nu\<{\A,\B}\.h = \nu\<{\A,\B}\.(i\<\A*a)\.f = a\.f = g.
\]

Let $h$ be a joint distribution on $\A\#\B$. We shall say that $a\:\A\>\B$ is a \NW{conditional}
IT for $h$ with respect to $\A$ whenever $h$ is generated by the marginal distribution
$\pi\<{\A,\B}\.h$ and the IT $a$, that is,
\[
  h = (i\<\A*a)\.\pi\<{\A,\B}\.h.
\]
Similarly, an IT $b\:\B\>\A$ such that
\[
  h = (b*i\<\B)\.\nu\<{\A,\B}\.h
\]
will be called a conditional IT for $h$ with respect to $\B$.

\subsection{Bayesian decision-making problems}
Suppose that there are fixed two
objects $\D$ and $\U$ in some category of ITs, namely, the object
of signals and the object of decisions, respectively. In a
decision-making problem with a prior distribution $f$ on $\D$ one
fixes some preorder $\BQ_f$ on the set of joint distributions on
$\D\#\U$ for which $\D$-marginal distribution coincides with $f$.

Informally, any joint distribution $h$ on $\D\#\U$ of this kind can be considered as a joint
distribution of a studied signal (with the distribution $f=\pi\<{\D,\U}\.h$ on $\D$) and a decision
(with the distribution $g=\nu\<{\D,\U}\.h$ on $\U$). The preorder $\BQ_f$ determines how good
is the ``correlation'' between studied signals and decisions.

Formally, an abstract \NW{decision problem with a prior information} is determined by a quadruple
$\TPL{\D,\U,f,\BQ_f}$, where $\D$ is an object of studied signals, $\U$ is an object of decisions
(or interpretations), $f\:\Z\>\D$ is a \NW{prior} distribution (or distribution \NW{a priori}), and
$\BQ_f$ is a preorder on the set of ITs $h\:\Z\>\D\#\U$ that satisfy the condition
$\pi\<{\D,\U}\.h=f$.

Furthermore, suppose that there is a fixed IT $a\:\D\>\R$ (which determines a measurement; $\R$
can be called an object of \NW{observations}). An IT $r\:\R\>\U$ is called \NW{optimal} (for the
IT $a$ with respect to $\BQ_f$) if the distribution $(i*r\.a)\.f$ is a maximal element with
respect to $\BQ_f$. The set of all optimal information transformers is denoted $\Opt_f(a\.f)$.

\begin{theorem} {\bf(Bayesian principle).} 
Let $f$ be a given prior distribution on $\D$, let $a\:\D\>\R$ be a fixed IT, and let $b\:\R\>\D$
be a conditional information transformer for $(i*a)\.f$ with respect to $\R$. Then the set of
optimal ITs $r\:\R\>\U$, namely, the set of optimal decision strategies for $f$ over $a\.f$
coincides with the set of optimal decision strategies for $b\.g$ over $g$, where $g=a\.f$:
\[
  \Opt_f(a\.f) = \Opt_{b\.g}(g).
\]
\end{theorem}

In a wide class of decision problems (e.g., in linear estimation
problems) an optimal IT $r$ happens to be deterministic and is
specified by the ``deterministic part'' of the IT $b$.

For many categories of information transformers (for example, stochastic, multivalued, and fuzzy
ITs{}~%
\ct{Borovkov, GolFil:MVMS, Gol:MSI, Gol:FuzDec}%
) an optimal decision strategy $r$ can be constructed ``pointwise'' according to the following
scheme. For the given ``result of observation'' $y\in\R$ consider the conditional (posterior)
distribution $b(y)$ for $f$ under a fixed $g=y$, and put
\[
  r(y) \Eqdef d_{b(y)},
\]
where $d_{b(y)}$ is an optimal decision with respect to the posterior distribution $b(y)$.

\subsection{Decision making problems in concrete categories of ITs}
\subsubsection{Stochastic ITs} 
\begingroup
\DefSth
Let us demonstrate here that the basic concepts of mathematical statistics are adequately described in
 terms of this IT-category. Namely, we shall verify that the concepts of distribution, conditional
distribution, etc. (introduced above in terms of IT-categories), in the category of stochastic ITs lead
to the corresponding classical concepts.

Indeed, any probability distribution $Q$ on a given measurable space $\A=\MS\A$ is uniquely
determined by the morphism $f\:\Z\>\A$ from the terminal object
$\Z=\TPL{\{0\},\;\bigl\{\EMPTYSET,\{0\}\bigr\}}$ (a one-point measurable space) such
that
\[
  P_f(0,A)=Q(A) \qquad \forall A\in\S_\A.
\]
In what follows we shall omit the first argument in $P_f(0,A)$ and write just $P_f(A)$ instead.

A statistical experiment is described by a family of probability measures $Q_\t$ on some measurable
space $\B$. This family is usually parametrized by elements of a certain set $\O_\A$. Sometimes
(especially when statistical problems with a prior information are studied) it is additionally assumed
that the set $\O_\A$ is equipped by some $\s$-algebra $\S_\A$ and that $Q_\t(B)$ is a
measurable function of $\t\in\O_\A$ for all $B\in\S_\B$ (and thus, $Q_\t(B)$ is a transition
probability function{}~%
\ct{Barra}).
Therefore, such statistical experiment is determined by the stochastic information transformer

$a\:\A\>\B$, where
\[
  P_a(\t,B)=Q_\t(B)\qquad \forall \t\in\O_\A, \quad \forall B\in\S_\B.
\]
In the case when no $\s$-algebra on the set $\O_\A$ is specified, one can put
$\S_\A=\spc{P}(\O_\A)$, that is, the $\s$-algebra of all the subsets of the set $\O_\A$. It is clear
that in this case the function $P_a(\t,B)=Q_\t(B)$ is a measurable function of $\t\in\O_\A$ for
every fixed $B\in\S_\B$ and thus (being a transition probability function), is described by a
stochastic IT $a\:\A\>\B$.

Note also, that any statistic, being a measurable function, is represented by a certain deterministic IT.
Decision strategies also correspond to deterministic ITs. At the same time, nondeterministic (mixed)
decision strategies are adequately represented by stochastic information transformers of general kind.

Now, let $f$ be some fixed distribution on $\A$ and let $a\:\A\>\B$ be some IT. The joint
distribution  $h$ on $\A\#\B$, generated by $f$ and $a$ (from the IT-categorical point of view, see
Section{}~\rs{Distr}) is
\[
  h = (i*a)\.f.
\]
It means that for every set $A\#B$, where $A\in\O_\A$ and $B\in\O_\B$,
\begin{eqnarray*}
  P_h(A\#B)
  &=&
  \Int_{\O_\A} P_{i*a}(\o,\,A\#B)\,P_f(d\o)
  \=
  \Int_{\O_\A} P_i(\o,A)\,P_a(\o,B)\,P_f(d\o)
  \=
  \Int_A P_a(\o,B)\,P_f(d\o).
\end{eqnarray*}
Thus we come to the well known classical expression for the generated joint distribution (see, for
example,{}~%
\ct{Barra}).

Now assume that $P_f$ is considered as some probability \NW{prior} distribution (or distribution
\NW{a priori}) on $\A$. Then for a given transition probability function $P_a$, a \NW{posterior}
(or \NW{conditional}) distribution $P_b(\o',\cdot)$ on $\A$ for a fixed $\o'\in\O_\B$ is
determined, accordingly to{}~%
\ct{Barra}
by a transition probability function $P_b(\o',A)$, $\o'\in\O_\B$, $A\in\S_\A$ such that
\[
\iftc
  P_h(A\#B) = \Int_B P_b(\o',A)\,P_g(d\o')\quad
  \forall A\in\S_\A, \forall B\in\S_\B,
\else
  P_h(A\#B) = \Int_B P_b(\o',A)\,P_g(d\o')\qquad
  \forall A\in\S_\A, \quad \forall B\in\S_\B,
\fi
\]
where
\[
  P_g(B) = \Int_{\O_\A} P_a(\o,B)\,P_f(d\o) \qquad \forall B\in\S_\B.
\]
It is easily verified that in terms of ITs the above expressions have the following forms:
\[
  h=(b*i)\.g,
\]
where
\[
  g=a\.f.
\]
This shows, that the classical concept of conditional distribution is adequately described by the
concept of conditional IT in terms of categories of information transformers.
\endgroup

\subsubsection{Linear stochastic ITs}
\begingroup
\DefLin
Note that in the category of linear information transformers
every IT is dominated (in the sense of the preorder relation
$\HP$) by a deterministic IT. Hence, according to
Proposition{}~\rp{DetMajor}, in any monotone decision-making
problem without loss of quality one can search optimal decision
strategies in the class of deterministic ITs.

According to section \rs{Distr} any joint distribution in
$\A\#\B$ is an IT $h\:\Z\>\A\#\B$, where $\Z$ is a terminal object
in the category of linear ITs, i.e., $\Z=\{0\}$ is a 0-dimensional
linear space. Thus, $h=\TPL{0,\S_{h}}$, where $\S_{h}$ is a
self-adjoint nonnegative operator in $h\:\Z\>\A\#\B$. Operator
$\S_{h}$ can be represented in the following ``matrix'' form:
\begin{equation} \lf{Joint:SLT}
 \S_{h} =
 \left(\begin{array}{cc}
   \S_{f} & \S_{f,g} \\ \S_{g,f} & \S_{g},
 \end{array}\right)
\end{equation}
where $\S_{f,g}$=$\S_{g,f}$.

It is shown in{}~%
\ct{Gol:RelInf}, that in the category of linear ITs for any joint
distribution there always exist conditional distributions.

\begin{theorem}
For any joint distribution $h\:\Z\>\A\#\B$ there exists conditional
information transformers $a\:\A\>\B$ and $b\:\B\>\A$.

Variants of conditional information transformers are given by the
formulas
\[
  a = \TPL{\S_{g,f}\S_f^-,\;\S_g - \S_{g,f}\S_f^-\S_{f,g}},
\]
\[
  b = \TPL{\S_{f,g}\S_g^-,\;\S_f - \S_{f,g}\S_g^-\S_{g,f}}.
\]
Here $A^-$ we denotes the pseudoinverse operator for
$A$~\ct{Pyt:PsInv}.

If $\S_f>0$ or $\S_g>0$ {\rm(}in this case this operator is
nonsingular and its pseudoinverse coincides with its inverse
$\S_f^{-}=\S_f^{-1}${\rm),} then the corresponding conditional
information transformer $a$ or $b$ is unique.
\end{theorem}

Thus in problems with a prior information one can apply Bayesian
principle. Its direct proof in the category of linear ITs as well
as the explicit expression for conditional information
transformers can
be found in{}~%
\ct{Gol:RelInf}.
\endgroup

\subsubsection{Multivalued ITs}
\begingroup
\DefMul
In the category of multivalued information transformers every IT is dominated (in the sense of the
partial order $\HP$) by a deterministic IT. Thus, in the monotone decision-making problem one can
search optimal decision strategies in the class of deterministic ones.

For every joint distribution in the category of multivalued ITs
there exist conditional
distributions{}~%
\ct{GolFil:MVMS}.

It is clear, that any joint distribution in $\A\#\B$, (i.e., an IT
$h\:\Z\>\A\#\B$) is specified by a subset $H$ of $\A\#\B$, since a
terminal object $\Z$ in the categories of multivalued ITs is a
1-element set. It is easy to see that for every joint distribution
in the category of multivalued ITs there exist conditional
distributions{}~%
\ct{GolFil:MVMS}.

\begin{theorem}%
For any joint distribution $H$ in $\A\#\B$, conditional
information transformers $a\:\A\>\B$ and $b\:\B\>\A$ always exist.
Some variants of conditional ITs are determined by the following
expressions:
\begin{eqnarray*}
  a x
  &=&
  \left\{\begin{array}{lll}
    {\displaystyle \SET{y\in\B|\TPL{x,y}\in H}}, &
    \mbox{if} & x\in p\<\A H,\\ [3mm]
    \B, & \mbox{if} & x\notin p\<\A H,
  \end{array}\right.
  \\[2mm]
  b y
  &=&
  \left\{\begin{array}{lll}
    {\displaystyle \SET{x\in\A|\TPL{x,y}\in H}}, &
    \mbox{if} & y\in p\<\B H,\\ [3mm]
    \A, & \mbox{if} & y\notin p\<\B H.
  \end{array}\right.
\end{eqnarray*}
Here $p\<\A H$ and $p\<\A H$ denote projections of $H$ on $\A$ and
$\B$ respectively, e.g., $p\<\A H\eqdef\SET{x\in\A|\exists
y\in\B\; \TPL{x,y}\in H}$.
\end{theorem}

Therefore, in decision problems with a prior information, the
Bayesian approach can be effectively applied.
\endgroup

\subsubsection{Categories of fuzzy information transformers}
\begingroup
\DefFuz

Here we define two categories of fuzzy information transformers
$\FMT$ and $\FPT$ that
correspond to different fuzzy theories{}~%
\ct{Gol:FuzDec}.

Objects of these categories are arbitrary sets and morphisms are
everywhere defined fuzzy maps, namely, maps that take an element
to a normed fuzzy set (a fuzzy set $A$ is normed if supremum of
its membership function $\x{A}$ is $1$). Thus, an information
transformer $a\:\A\>\B$ is defined by a membership function
$\x{ax}(y)$ which is interpreted as the grade of membership of an
element $y\in\B$ to a fuzzy set $ax$ for every element $x\in\A$.

{\it The category $\FMT$}.
Suppose $a\:\A\>\B$ and $b\:\B\>\C$ are some fuzzy maps. We define
their \NW{composition} $b\.a$ as follows: for every element
$x\in\A$ put
\[
  \x{(b\.a)x}(z) \Eqdef \sup_{y\in\B}\min\(\x{ax}(y),\;\x{by}(z)\).
\]

For a pair of fuzzy information transformers $a\:\D\>\A$ and
$b\:\D\>\B$ with the common source $\D$, we define their
\NW{product} as the IT that acts from $\D$ to the Cartesian
product $\A\#\B$, such that
\[
  \x{(a*b)x}(y,z) \Eqdef \min\(\x{ax}(y),\;\x{by}(z)\).
\]

{\it The category $\FPT$}.
Define the \NW{composition} and the \NW{product} by the following
expressions:
\[
  \x{(b\.a)x}(z) \Eqdef \sup_{y\in\B}\(\x{ax}(y)\;\x{by}(z)\),
\]
\[
  \x{(a*b)x}(y,z) \Eqdef \x{ax}(y)\;\x{by}(z).
\]

In the both defined above categories of fuzzy information
transformers the subcategory of \NW{deterministic} ITs is
(isomorphic to) the category of sets $\Set$. Let $g\:\A\>\B$ be
some map (morphism in $\Set$). Define the corresponding fuzzy IT
(namely, a fuzzy map, which is obviously, everywhere defined)
$\tilde{g}\:\A\>\B$ in the following way:
\[
  \x{\tilde{g}(x)}(y)\Eqdef\d{g(x),y}
  =
  \CASE 1 & g(x)=y; 0 & g(x)\neq y.;
\]

Concerning the choice of accuracy relation, note, that in these
IT\_categories, like in the category of multivalued ITs, apart
from the trivial accuracy relation one can put for $a,b\:\A\>\B$
\[
  a \HP b \Iffdef \forall x\in\A\q \forall y\in\B\q \x{ax}(y)\LE \x{bx}(y).
\]
In each fuzzy IT-category these two choices lead to two different informativeness relations, namely
the strong and the weak ones.

Like in the categories of linear and multivalued ITs discussed above, monotone decision-making
problems admit restriction of the class of optimal decision strategies to deterministic ITs without loss
of quality.

It is clear, that any joint distribution in $\A\#\B$, (i.e., an IT
$h\:\Z\>\A\#\B$) is, in fact a normed fuzzy subset of $\A\#\B$,
since a terminal object $\Z$ in the categories of fuzzy ITs is a
1-element set. Denote this fuzzy set by $H$. It is shown in{}~%
\ct{Gol:FuzDec, Gol:MSI} that for every joint distribution in the
categories of fuzzy ITs there exist conditional distributions.

\begin{theorem}%
For any joint distribution $H$ in $\A\#\B$, conditional
information transformers $a\:\A\>\B$ and $b\:\B\>\A$ always exist.
Some variants of conditional ITs are determined by the following
expressions:

{\rm(a)} in the category $\FMT$
\begin{eqnarray*}
  \x{a x}(y)&=&\x H(x,y),
  \\[2mm]
  \x{b y}(x)&=&\x H(x,y);
\end{eqnarray*}
\begin{eqnarray*}
  \x{a x}(y)
  &=&
  \left\{\begin{array}{lll}
    \x H(x,y), &
    \mbox{if} & \x H(x,y) \neq \sup\limits_z\x H(x,z),\\ [3mm]
    1, & \mbox{if} & \x H(x,y) = \sup\limits_z\x H(x,z),
  \end{array}\right.
  \\[2mm]
  \x{b y}(x)
  &=&
  \left\{\begin{array}{lll}
    \x H(x,y), &
    \mbox{if} & \x H(x,y) \neq \sup\limits_z\x H(z,y), \\ [3mm]
    1, & \mbox{if} & \x H(x,y) \neq \sup\limits_z\x H(z,y).
  \end{array}\right.
\end{eqnarray*}

{\rm(b)} in the category $\FPT$
\begin{eqnarray*}
  \x{a x}(y)
  &=&
  \left\{\begin{array}{lll}
    {\displaystyle\frac{\x H(x,y)}{\sup\limits_z\x H(x,z)}}, &
    \mbox{if} & \sup\limits_z\x H(x,z)\not=0,\\ [3mm]
    1, & \mbox{if} & \sup\limits_z\x H(x,z)=0,
  \end{array}\right.
  \\[2mm]
  \x{b y}(x)
  &=&
  \left\{\begin{array}{lll}
    {\displaystyle\frac{\x H(x,y)}{\sup\limits_z\x H(z,y)}}, &
    \mbox{if} & \sup\limits_z\x H(z,y)\not=0, \\ [3mm]
    1, & \mbox{if} & \sup\limits_z\x H(z,y)=0.
  \end{array}\right.
\end{eqnarray*}
\end{theorem}

This allows Bayesian approach and makes use of Bayesian principle in
decision problems with a prior
information for fuzzy ITs{}~%
\ct{Gol:FuzDec, Gol:MSI, Gol:FuzLog, Gol:MSIII},
where connections between fuzzy decision problems and the underlying fuzzy logic are studied.

The present paper, while very much a first step, lays the basis for number of further applications. In paper \ct{Gol:MSIII} we propose some realizations for above categories, which we belive can be the basis for some interesting new directions in quantum computation and bioinformatics.

\endgroup



\newcommand{\NUMBER}{N }    

\makeatletter
\def\@biblabel#1{#1.}
\makeatother

\def\SETLANG#1{%
  \if#1R%
    \def\VOL{'}\def\PAG{'}%
  \else%
    \def\VOL{vol}\def\PAG{P}%
  \fi%
}

\def\ARTX A=#1 T=#2 J=#3 Y=#4 V=#5 N=#6 P=#7 L=#8{%
#1 #2. {\it#3,\/} #4,
{\if?#5\else {\bf#5,}\fi}
{\if?#6\else N~#6,\fi}
{\if?#7\else P.~#7.\fi}}

\newcommand{\Bibitem}[1]{\bibitem{#1}}

\def\ARTE#1 A=#2 T=#3 J=#4 Y=#5 V=#6 N=#7 P=#8 \par{\Bibitem{#1}
\ARTX A=#2 T=#3 J=#4 Y=#5 V=#6 N=#7 P=#8 L=E}

\def\ART#1 A=#2 T=#3 J=#4 Y=#5 V=#6 N=#7 P=#8 L=#9{\Bibitem{#1}
\ARTX A=#2 T=#3 J=#4 Y=#5 V=#6 N=#7 P=#8 L=#9}


\def\BOKX A=#1 T=#2 J=#3 Y=#4 =#5 \par{#1 {\it#2.\/} #3, #4#5 \par}

\def\BOK#1 A=#2 T=#3 J=#4 Y=#5 =#6 \par{\Bibitem{#1}
  \BOKX A=#2 T=#3 J=#4 Y=#5 =#6. \par}

\def\MET#1 A=#2 T=#3 J=#4 H=#5 Y=#6 P=#7 \par{\Bibitem{#1}
#2 #3. {\it#4,\/} #5, #6, pp.~#7.}



\begin{thebibliography}{99}



\ARTE{Zad}
A=Zadeh~L. A.  
T=The concept of a linguistic variable and its 
application to approximate reasoning. Parts 1, 2, and 3
J=Inf. Sci. 
Y=1975 V=8  N=? P=199-259, 301--357; 1975, {\bf 9}, P.~43--80


\BOK{Dub:FSI}  
A=Dubois D., Prade H. 
T=Fuzzy Sets and Systems: Theory and 
Applications. 
J=New York: Academic Press 
Y=1979  =

\ARTE{Dub:FSII}
A=Dubois D., Prade H. 
T=Fuzzy sets in approximate reasoning, Part I: 
Inference with possibility distributions
J=Fuzzy Sets Syst.
Y=1991 V=40 N=? P=143--202

\ARTE{Gol:FuzDec}
 A=Golubtsov~P.V.
 T=Theory of Fuzzy Sets as a Theory of Uncertainty and
   Decision-Making Problems in Fuzzy Experiments
 J=Probl. Inf. Transm.
 Y=1994 V=30 N=3 P=232--250


\ARTE{Gol:MSI}
 A=Golubtsov~P.V., Moskaliuk S.S.
 T=Bayesian Decisions and Fuzzy Logic
 J=Preprint ESI Num. 626 Vienna, Austria
  Y=1998 V=? N=? P=1--8

\ARTE{Manes:ClassFuzzy}
 A=Manes~E.G.
 T=A Class of Fuzzy Theories
 J=J. Math. Anal. and Appl. Y=1982 V=85 N=? P=409--451

\BOK{Black:Gir}                  
 A=Blackwell~D., Girschick~M.A.
 T=Theory of Games and Statistical Decisions
 J=New York: Wiley Sons Y=1954 =


\BOK{Barra}
 A=Barra~J.-R.
 T=Notions Fondamentales de Statistique
   Math\'ematique
 J=Paris: Dunod
 Y=1974 = 




\ARTE{Sackst:StatEq}
 A=Sacksteder~R.
 T=A Note on Statistical Equivalence
 J=Ann. Math. Statist.
 Y=1967 V=38 N=3 P=787--795

\ARTE{MorseSackst:StatIso}
 A=Morse~N., Sacksteder~R.
 T=Statistical Isomorphism
 J=Ann. Math. Statist. Y=1966 V=37 N=2 P=203--214

\BOK{Chentsov}
 A=Chentsov~N.N.
 T=Statistical Decision Rules and Optimal Inference
   {\rm [in Russian]}
 J=Moscow: Nauka Y=1972 =    

\ARTE{Chentsov:CatMatStat}
 A=Chentsov~N.N.
 T=Categories of Mathematical Statistics
 J=Dokl. Akad. Nauk SSSR
 Y=1965 V=164 N=3 P=511--514

\ARTE{Gol:MSAlg}
 A=Golubtsov~P.V.
 T=Measurement Systems: Algebraic
   Properties and Informativity
 J=Pattern Recognition and Image Analysis
 Y=1991 V=1 N=1 P=77--86

\BOK{MacLane}
 A=MacLane~S.
 T=Categories for the Working Mathematician
 J=New York: Springer Y=1971 =

\BOK{Herrlich:Strecker}
 A=Herrlich~H., Strecker~G.E.
 T=Category Theory
 J=Boston: Allyn and Bacon Y=1973 =

\BOK{Arrows}
 A=Arbib~M.A., Manes~E.G.
 T=Arrows, Structures and Functors
 J=New York: Academic Press Y=1975 =

\BOK{Goldblatt}
 A=Goldblatt~R.
 T=Topoi. The Categorial Analysis of Logic
 J=Amsterdam: North-Holland Y=1979 =

\ARTE{Gol:InfCatLinSys}
  A=Golubtsov~P.V.
  T=Informativity in the Category of Linear Measurement Systems
  J=Probl. Inf. Transm.
  Y=1992 V=28 N=2 P=125--140

\ARTE{Gol:RelInf}
  A=Golubtsov~P.V.
  T=Relative Informativity and a Priori Information in the Category
    of Linear Information Transformers
  J=Probl. Inf. Transm.
  Y=1995 V=31 N=3 P=195--215

\ARTE{GolFil:MVMS}
 A=Golubtsov~P.V., Filatova~S.A.
 T=Multivalued Mea\-su\-re\-ment-Computer Systems
 J=Mat. Model. Y=1992 V=4 N=7 P=79--94

\ARTE{Gol:MV_Inf}
  A=Golubtsov~P.V.
  T=Informativity in the Category of Multivalued Information Transformers
  J=Probl. Inf. Transm.
  Y=1998 V=34 N=3 P=259--276

\MET{Gol:FuzLog}
 A=Golubtsov~P.V.
 T=Fuzzy Logical Semantics of Bayesian Decision Making
 J=Applications of Fuzzy Logic Technology II. SPIE Proceedings. Vol.~2493.
 H=Orlando, Florida Y=1995 P=228--239

\MET{Gol:CatInfTran}
 A=Golubtsov~P.V.
 T=Categories of Information Transformers and the Concept of Informativity
 J=Proc. Int. Conf. on Informatics and Control (ICI\&C'97). Vol.~2
 H=St.\_Petersburg, Russia Y=1997. P=512--517

\ARTE{Gol:Ax-IT}
  A=Golubtsov~P.V.
  T=Axiomatic Description of Categories of Information
    Transformers
  J=Probl. Inf. Transm.
  Y=1999 V=35 N=3 P=80--99


\ARTE{Gol:MS}
  A=Golubtsov~P.V., Moskaliuk~S.S.
  T= Method of Additional Structures on the Objects of a Monoidal 
  Kleisli Category as a Background for Information Transformers
  J=Hadronic Journal  
  Y=2002 V=25 N=2 P=179--238



\ARTE{Kleisli}
 A=Kleisli~H.
 T=Every Standard Construction is
   Induced by a Pair of Adjoint Functors
 J=Proc. Amer. Math. Soc. Y=1965 V=16 N=? P=544--546

\BOK{BarrWells:TTT}
 A=Barr M., Wells C.
 T=Toposes, Triples and Theories
 J=New York: Springer Y=1984 =

\BOK{Borovkov}
 A=Borovkov~A.A.
 T=Mathematical Statistics. Supplementary Chapters. {\rm[in Russian]}
 J=Moscow: Nauka
 Y=1984 = 

\MET{Black:Comp}
 A=Blackwell~D.
 T=Comparison of Experiments
 J=Proc. Second Berkeley Sympos. on Mathematical Statistics and
   Probability.
 H=University of California Press Y=1951 P=93--102

\ARTE{Black:EquComp}
 A=Blackwell~D.
 T=Equivalent Comparison of Experiments
 J=Ann. Math. Statist.
 Y=1953 V=24 N=2 P=265--272


\MET{Giry}
 A=Giry~M.
 T=A Categorical Approach to Probability Theory
 J=Categorical Aspects of Topology and Analysis.
   Lecture Notes in Mathematics. \NUMBER~915.
 H=Berlin: Springer-Verlag Y=1982 P=68--85

\BOK{Neveu}
 A=Neveu~J.
 T=Bases Math\'ematiques du Calculus des Probabilit\'es
 J=Paris: Masson et Cie Y=1964 =


\ARTE{Pyt:ZadRed}
 A=Pyt'ev~Yu.P.
 T=Reduction Problems in Experimental Research
 J=Mat. Sb. Y=1983 V=120 N=2 P=240--272

\ARTE{Pyt:PsInv}
 A=Pyt'ev~Yu.P.
 T=Pseudoinverse operator. Properties and applications
 J=Mat. Sb. Y=1982 V=118 N=5 P=19--49

\ARTE{Goguen:L-Fuzzy}
 A=Goguen~J.A.
 T=L-Fuzzy Sets
 J=J. Math. Anal. and Appl.
 Y=1967 V=18 N=? P=145--174

\ARTE{Goguen:V-Sets}
 A=Goguen~J.A.
 T=Categories of V-Sets
 J=Bull. Amer. Math. Soc.
 Y=1969 V=75 N=? P=622--624

\ARTE{Goguen:Concept}
 A=Goguen~J.A.
 T=Concept Representation in Natural and Artificial
   Languages: Axioms, Extensions, and Applications
   for Fuzzy Sets
 J=Int. J. Man-Machine Studies Y=1974 V=6 N=? P=513--561

\ARTE{FuzzyMach}
 A=Arbib~M.A., Manes~E.G.
 T=Fuzzy Machines in a Category
 J=Bull. Austral. Math. Soc. Y=1975 V=13 N=? P=169--210

\MET{FuzzyMorph}
 A=Arbib~M.A., Manes~E.G.
 T=Fuzzy Morphisms in Automata Theory
 J=Category Theory Applied to Computation and Control.
   Lecture Notes in Computer Science. No.~25.
 H=New York: Springer-Verlag
 Y=1976 P=80--86

\BOK{Manes:AlgTh}
 A=Manes~E.G.
 T=Algebraic Theories
 J=New York: Springer Y=1976 =


\ARTE{Gol:MSIII}
  A=Moskaliuk S.S., Demchynska M. I. 
  T=Information Transformers in Bioinformatics
  J=Hadronic Journal  
  Y= V= N= P= (to appear) 



\end{thebibliography}
\end{document}